\documentclass[12pt,letterpaper]{amsart}

\newcounter{myenum}

\newcommand{\indentalign}{\hspace{0.3in}&\hspace{-0.3in}}
\newcommand{\la}{\langle}
\newcommand{\ra}{\rangle}

\newcommand{\I}{\operatorname{Im}}

\newcommand{\nls}{\textsc{NLS}}
\newcommand{\defeq}{\stackrel{\rm{def}}{=}}

\usepackage{amssymb}
\usepackage{amsthm}
\usepackage{amsxtra}

\newtheorem{theorem}{Theorem}
\newtheorem{definition}[theorem]{Definition}
\newtheorem{proposition}[theorem]{Proposition}
\newtheorem{lemma}[theorem]{Lemma}
\newtheorem{corollary}[theorem]{Corollary}

\theoremstyle{remark}
\newtheorem{remark}[theorem]{Remark}

\newcommand{\cR}{\mathbb{R}}

\newcommand{\ep}{\epsilon}

\newcommand{\lap}{\triangle}
\newcommand{\grad}{\nabla}

\numberwithin{equation}{section}
\numberwithin{theorem}{section}

\title[Scattering of 3d cubic NLS]
{A sharp condition for scattering of the radial 3d cubic nonlinear Schr\"odinger equation}
\author{Justin Holmer}
\address{University of California, Berkeley}
\author{Svetlana Roudenko}
\address{Arizona State University}

\begin{document}

\maketitle

\begin{abstract}
We consider the problem of identifying sharp criteria under which
radial $H^1$ (finite energy) solutions to the focusing 3d cubic
nonlinear Schr\"odinger equation (NLS) $i\partial_t u + \Delta u +
|u|^2u=0$ scatter, i.e.\ approach the solution to a linear
Schr\"odinger equation as $t\to \pm \infty$.   The criteria is
expressed in terms of the scale-invariant quantities
$\|u_0\|_{L^2}\|\nabla u_0\|_{L^2}$ and $M[u]E[u]$, where $u_0$
denotes the initial data, and $M[u]$ and $E[u]$ denote the
(conserved in time) mass and energy of the corresponding solution
$u(t)$.  The focusing NLS possesses a soliton solution $e^{it}Q(x)$,
where $Q$ is the ground-state solution to a nonlinear elliptic
equation, and we prove that if $M[u]E[u]<M[Q]E[Q]$ and
$\|u_0\|_{L^2}\|\nabla u_0\|_{L^2} < \|Q\|_{L^2}\|\nabla Q\|_{L^2}$,
then the solution $u(t)$ is globally well-posed and scatters.  This
condition is sharp in the sense that the soliton solution
$e^{it}Q(x)$, for which equality in these conditions is obtained, is
global but does not scatter.  We further show that if
$M[u]E[u]<M[Q]E[Q]$ and $\|u_0\|_{L^2}\|\nabla u_0\|_{L^2} >
\|Q\|_{L^2}\|\nabla Q\|_{L^2}$, then the solution blows-up in finite
time.  The technique employed is parallel to that employed by
Kenig-Merle \cite{KM06a} in their study of the energy-critical NLS.
\end{abstract}

\section{Introduction}

Consider the cubic focusing nonlinear Schr\"odinger (NLS) equation on $\mathbb{R}^3$:
\begin{equation}
 \label{E:NLS}
i\partial_t u +\Delta u + |u|^2u=0,
\end{equation}
where $u=u(x,t)$ is complex-valued and $(x,t)\in \mathbb{R}^3\times
\mathbb{R}$. The initial-value problem posed with initial-data
$u(x,0)=u_0(x)$ is locally well-posed in $H^1$ (see Ginibre-Velo
\cite{GV79}; standard reference texts are Cazenave \cite{Caz-book},
Linares-Ponce \cite{LP-book}, and Tao \cite{Tao-book}). Such
solutions, during their lifespan $[0,T^*)$ (where $T^*=+\infty$ or
$T^*<+\infty$), satisfy mass conservation $M[u](t)=M[u_0]$, where
$$
M[u](t) = \int |u(x,t)|^2 \,dx,
$$
and energy conservation $E[u](t)=E[u_0]$, where
$$
E[u](t) = \frac12\int |\nabla u(x,t)|^2 \, dx - \frac14\int
|u(x,t)|^4 \,dx
$$
(and we thus henceforth denote these quantities $M[u]$ and $E[u]$
respectively, with no reference to the time $t$).

The equation also has several invariances, among them (in each of
the following cases, $\tilde u$ is a solution to \eqref{E:NLS} if
and only if $u$ is a solution to \eqref{E:NLS}):
\begin{itemize}
\item \textit{Spatial translation}.  For a fixed $x_0\in \mathbb{R}^3$, let $\tilde u(x,t) = u(x+x_0,t)$.
\item \textit{Scaling}.  For a fixed $\lambda \in (0,+\infty)$, let $\tilde u(x,t) = \lambda u(\lambda x, \lambda^2 t)$.
\item \textit{Galilean phase shift}.  For a fixed $\xi_0\in \mathbb{R}^3$, let $\tilde u(x,t) = e^{ix\xi_0}e^{-it\xi_0^2}u(x-2\xi_0t,t)$.
\end{itemize}
The scale-invariant Sobolev norm is $\dot H^{1/2}$, although we find
it more useful, as described below, to focus on the scale invariant
quantities $\|u(t)\|_{L^2}\|\nabla u(t)\|_{L^2}$ and $M[u]E[u]$. The
Galilean invariance leaves only the $L^2$ norm invariant, while
translation leaves all Sobolev norms invariant.  We note that these
two symmetries do not preserve radiality, while the scaling symmetry
does.

The nonlinear elliptic equation
\begin{equation}
\label{E:Q} -Q + \Delta Q + |Q|^2Q=0, \qquad Q=Q(x), \qquad x\in
\mathbb{R}^3,
\end{equation}
has an infinite number of solutions in $H^1$. Among these there is
exactly one solution of minimal mass\footnote{In view of the connection
between solutions $Q$ to \eqref{E:Q} and solutions $u(t)=e^{it}Q$ to
\eqref{E:NLS}, and the fact that $\|u(t)\|_{L^2}\|\nabla
u(t)\|_{L^2}$ is a scale invariant quantity for solutions $u(t)$ to
\eqref{E:NLS}, it might be more natural to classify the family of
solutions $Q$ to \eqref{E:Q} in terms of the quantity
$\|Q\|_{L^2}\|\nabla Q\|_{L^2}$ rather than the mass.  
However, any solution $Q$ to
\eqref{E:Q} must satisfy the Pohozhaev identity $\|Q\|_{L^2}\|\nabla
Q\|_{L^2} = \sqrt{3}\|Q\|_{L^2}^2$, and thus the two 
classifications are equivalent.}, called the
\textit{ground-state} solution, and it is positive (real-valued),
radial, smooth, and exponentially decaying (see Appendix B of Tao's
text \cite{Tao-book} for exposition).  We henceforth denote by $Q$
this ground-state solution.  If we let $u(x,t) = e^{it}Q(x)$, then
$u$ is a solution to \eqref{E:NLS}, and is called the
\textit{standard soliton}.  A whole family of
soliton solutions to \eqref{E:NLS} can be built from the standard
soliton via the invariances of the NLS equation \eqref{E:NLS}:
\begin{equation}
\label{E:solitons} u(x,t) = e^{it}e^{ix\cdot
\xi_0}e^{-it|\xi_0|^2}\lambda \,
u(\lambda(x-(x_0+2\xi_0t)),\lambda^2 t).
\end{equation}
The standard soliton has the property that the quantities
$\|u_0\|_{L^2}\|\nabla u_0\|_{L^2}$ and $M[u]E[u]$ are minimal among
all solitons \eqref{E:solitons}. Indeed, these quantities are
independent of translation and scaling, and the introduction of a
Galilean phase shift only increases their values. Since  solutions
to the linear Schr\"odinger equation completely disperse (spread
out, and shrink in a variety of spatial norms) as $t\to \pm \infty$,
the soliton solutions  by their definition do not scatter (approach
a solution of the linear Schr\"odinger equation). Indeed, soliton
solutions represent a perfect balance between the focusing forces of
the nonlinearity and the dispersive forces of the linear component.

The basic line of thought in the subject, motivated by heuristics
(Soffer \cite{Soffer}), rigorous partial results (Tao \cite{T04,T06}),
numerical simulation (Sulem-Sulem \cite{SS99}), and analogy with the
completely integrable one-dimensional case, is that a solution of
\eqref{E:NLS} either completely disperses as $t\to \infty$ (linear
effects dominate), blows-up in finite time (nonlinear effects
dominate) or the solution resolves into a sum of solitons
propagating in different directions or at different speeds plus
dispersive radiation as $t \to \infty$ (nonlinear effects and linear
effects balance).  Since the smallest value of
$\|u_0\|_{L^2}\|\nabla u_0\|_{L^2}$ among all soliton solutions is
$\|Q\|_{L^2}\|\nabla Q\|_{L^2}$, it seems reasonable to conjecture,
even for nonradial data, that if $\|u_0\|_{L^2}\|\nabla
u_0\|_{L^2}<\|Q\|_{L^2}\|\nabla Q\|_{L^2}$, then the solution
scatters provided we can rule out blow-up.  Ruling out blow-up in
this situation is straightforward provided $M[u]E[u]<M[Q]E[Q]$ using
the conservation of mass and energy and a result of M.\ Weinstein
stating that an appropriate Gagliardo-Nirenberg inequality is
optimized at $Q$.  The main result of this paper is the resolution
of this conjecture under the assumption of radial data, which
appears below as Theorem \ref{T:main}\eqref{I:global}(b).

\begin{theorem}
\label{T:main} Let $u_0\in H^1$ be radial and let $u$ be the
corresponding solution to \eqref{E:NLS} in $H^1$ with maximal
forward time interval of existence $[0,T)$. Suppose
$M[u]E[u]<M[Q]E[Q]$.
\begin{enumerate}
\item
 \label{I:global}
If $\|u_0\|_{L^2}\|\nabla u_0\|_{L^2}<\|Q\|_{L^2}\|\nabla
Q\|_{L^2}$, then
\begin{enumerate}
\item
$T=+\infty$ (the solution is globally well-posed in $H^1$), and
\item
$u$ scatters in $H^1$. This means that there exists $\phi_+ \in H^1$
such that
$$
\lim_{t\to +\infty} \|u(t)-e^{it\Delta}\phi_+\|_{H^1} =0\, .
$$
\end{enumerate}
\item
 \label{I:blow-up}
If $\|u_0\|_{L^2}\|\nabla u_0\|_{L^2}>\|Q\|_{L^2}\|\nabla
Q\|_{L^2}$, then $T<+\infty$ (the solution blows-up in finite time).
\end{enumerate}
\end{theorem}
It is straightforward to establish as a corollary the same result
for negative times: take the complex conjugate of the equation and
replace $t$ by $-t$.  Since the hypotheses in Theorem \ref{T:main}
\eqref{I:global}\eqref{I:blow-up} apply to $u_0$ if and only if they
apply to $\bar u_0$, we obtain that the hypotheses of
\eqref{I:global} imply that $u$ scatters both as $t\to +\infty$ and
$t\to -\infty$ and the hypotheses of \eqref{I:blow-up} imply that
$u$ blows-up both in finite positive time and in finite negative
time.  An interesting open question is whether or not there exist
solutions $u$ with $M[u]E[u]\geq M[Q]E[Q]$ that exhibit different
behavior in the positive and negative directions.

The proof of Theorem \ref{T:main}\eqref{I:global}(b) is based upon
ideas in Kenig-Merle \cite{KM06a}, who proved an analogous statement
for the energy-critical NLS.

The key dynamical quantity in the proof of Theorem \ref{T:main} is a
localized variance $\| x u(t)\|_{L^2({|x|\leq R})}$. The
\textit{virial identity} states that if $\|xu_0\|_{L^2}<\infty$,
then $u$ satisfies
\begin{equation}
 \label{E:virial}
\partial_t^2 \int |x|^2 |u(x,t)|^2 \, dx = 24E[u] - 4\|\nabla
u(t)\|_{L_x^2}^2.
\end{equation}
We use a localized version of this identity in both the proof of
Theorem \ref{T:main}\eqref{I:blow-up} and the rigidity lemma (see \S
\ref{S:rigidity}) giving Theorem \ref{T:main}\eqref{I:global}(b). On
a heuristic level (keeping in mind that $\|u(t)\|_{L^2}$ is
conserved), under the hypotheses of Theorem
\ref{T:main}\eqref{I:global}, the right side of \eqref{E:virial} is
strictly positive, which pushes the variance $\|xu(t)\|_{L^2}$ to
$+\infty$ as $t\to +\infty$, which says roughly that the mass of $u$
is being redistributed to large radii, meaning that it
``disperses'', and we expect the effect of the nonlinearity to
diminish and scattering to occur.  On the other hand, under the
hypotheses of Theorem \ref{T:main}\eqref{I:blow-up}, the right side
of \eqref{E:virial} is strictly negative, which pushes the variance
$\|xu(t)\|_{L^2}$ to $0$ in finite time, meaning that all the mass
of $u$ concentrates at the origin and we expect blow-up.  We do not
use \eqref{E:virial} directly, however, for two reasons.  First, it
requires the additional hypothesis that the initial data has finite
variance--an assumption we would like to avoid.  Secondly, in the
case of the scattering argument, we don't see a method for
\textit{proving} scattering given only the strict convexity (in
time) of the variance and its divergence to $+\infty$, although it
is \textit{heuristically} consistent with scattering.  The problem
is that large variance can be produced by a very small amount of
mass moving to very, very large radii, while still leaving a
significant amount of mass at small radii.  Therefore, to prove the
scattering claim in Theorem \ref{T:main}\eqref{I:global}(b), we
instead use a localized virial identity, as Kenig-Merle \cite{KM06a}
did, involving a localized variance.  If a very small amount of mass
moved to very, very large radii, it would not affect the
\textit{localized} variance dramatically.

For the 3d cubic defocusing NLS
$$
i\partial_t u + \Delta u -|u|^2u=0
$$
scattering has been established for all $H^1$ solutions (regardless
of ``size'') even for nonradial data by Ginibre-Velo \cite{GV85}
using a Morawetz inequality. This proof was simplified by
Colliander-Keel-Staffilani-Takaoka-Tao \cite{CKSTT04} using a new
interaction Morawetz inequality they discovered. These Morawetz
estimates, however, are not positive definite for solutions to the
focusing equation \eqref{E:NLS}, and thus, cannot be applied
directly to our problem.  It remains open whether or not one could
prove suitable bounds on the nonpositive terms to recover the
results of this paper.

For \eqref{E:NLS}, Tao \cite{T04} proved a few results in the
direction of the soliton resolution conjecture, assuming the
solution is radial and global (has globally bounded $H^1$ norm). It
is shown that for large data, radial solutions asymptotically split
into (i) a (smooth) function localized near the origin (which is
either zero or has a non-zero mass and energy and obeys an asymptotic
Pohozhaev identity), (ii) a radiation term evolving by the linear
Schr\"odinger flow, and (iii) an error term (approaching zero in the
$\dot{H}^1$ norm). Further results for mass supercritical, energy
subcritical NLS equations in higher dimensions ($N \geq 5$) were
established by Tao in \cite{T06}.

The equation \eqref{E:NLS} frequently arises, often in more complex
forms, as a model equation in physics.  In 2d, it appears as a model
in nonlinear optics -- see Fibich \cite{F} for a review.  When
coupled with a nonlinear wave equation, it arises as the Zakharov
system \cite{27} in plasma physics.  According to \cite{26} p.7, in
the mass supercritical case ``the most important partial case $p=3$,
$d=3$ corresponds to the subsonic collapse of Langmuir waves in
plasma''.  Furthermore, \eqref{E:NLS} arises as a model for the
Bose-Einstein condensate (BEC) in condensed matter physics.  There,
it appears as the Gross-Pitaevskii (GP) equation (in 1d, 2d, and
3d), which is \eqref{E:NLS} with a (real) potential $V=V(x)$:
\begin{equation}
\label{E:GP}
i\partial_t u + \Delta u -Vu + a|u|^2u=0 \; .
\end{equation}
It is derived by mean-field theory approximation (see Schlein
\cite{S}), and $|u(x,t)|^2$ represents the density of the condensate
at time $t$ and position $x$.  The coefficient $a$ in the
nonlinearity is governed by a quantity called the $s$-scattering
length.  Some elements used in recent experiments (${}^7$Li,
${}^{85}$Rb, ${}^{133}$Cs) posses a negative $s$-scattering length
in the ground state and are modeled by \eqref{E:GP} with $a<0$.
$V(x)$ is an external trapping potential imposed by a system of
laser beams and is typically taken to be harmonic $V(x)=\beta
|x|^2$.  These ``unstable BECs'' (where $a<0$) have been
investigated experimentally recently (see the JILA experiments
\cite{JILA}) and a number of theoretical predictions have been
confirmed, including the observation of ``collapse events''
(corresponding to blow-up of solutions to \eqref{E:GP}).   A few
articles have appeared (for example \cite{BAK}) in the physics
literature discussing the critical number of atoms required to
initiate collapse.  The ``critical number of atoms'' corresponds to
``threshold mass $M[u]$'' in our terminology, and connects well with
the mathematical investigations in this paper.

The format of this paper is as follows. In \S\ref{S:localtheory}, we
give a review of the Strichartz estimates, the small data theory,
and the long-time perturbation theory.  We review properties of the
ground state profile $Q$ in \S\ref{S:groundstate} and recall its
connection to the sharp Gagliardo-Nirenberg estimate of M.\
Weinstein \cite{W83}.  In \S\ref{S:dichotomy}, we introduce the
local virial identity and prove Theorem \ref{T:main} except for the
scattering claim in part \eqref{I:global}(b). In
\S\ref{S:compactness}-\ref{S:rigidity}, we prove Theorem
\ref{T:main}\eqref{I:global}(b). This is done in two stages,
assuming that the threshold for scattering is strictly below the one
claimed. First, in \S\ref{S:compactness}, we construct a solution
$u_\text{c}$ (a ``critical element'') that stands exactly at the
boundary between scattering and nonscattering. This is done using a
profile decomposition lemma in $\dot H^{1/2}$, obtained by extending
the $\dot H^1$ methods of Keraani \cite{K01}. We then show that time
slices of $u_\text{c}(t)$, as a collection of functions in $H^1$,
form a precompact set in $H^1$ (and thus $u_\text{c}$ has something
in common with the soliton $e^{it}Q(x)$). This enables us to prove
that $u_\text{c}$ remains localized \textit{uniformly in time}. In
\S\ref{S:rigidity}, this localization is shown to give a strict
convexity (in time) of a localized variance which leads to a
contradiction with the conservation of mass at large times. In
\S\ref{S:extensions}, we explain how Theorem 4.2 should
carry over to more general nonlinearities and general dimensions
(mass supercritical and energy subcritical cases) of NLS equations.

{\bf Acknowledgement.} J.H. is partially supported by an
NSF postdoctoral fellowship.
S.R. would like to thank Mary and Frosty Waitz for their great
hospitality during her visits to Berkeley. We both thank Guixiang Xu
for pointing out a few misprints and the referee for helpful
suggestions.

\section{Local theory and Strichartz estimates}
\label{S:localtheory}

We begin by recalling the relevant Strichartz estimates (e.g., see
Cazenave \cite{Caz-book},  Keel-Tao \cite{KT98}). We say that
$(q,r)$ is $\dot H^s$ Strichartz admissible (in 3d) if
$$
\frac2q+\frac3r=\frac32-s.
$$
Let
$$
\|u\|_{S(L^2)} = \sup_{\substack{(q,r)\; L^2 \text{ admissible} \\
2\leq r \leq 6, \;  2\leq q \leq \infty}} \|u\|_{L_t^qL_x^r}.
$$
In particular, we are interested in $(q,r)$ equal to
$(\frac{10}{3},\frac{10}3)$ and $(\infty, 2)$.  Define\footnote{For
some inequalities, the range of valid exponents $(q,r)$ can be
extended. The Kato inequality \eqref{E:Kato} imposes the most
restrictive assumptions that we incorporate into our definitions of
$S(\dot H^{1/2})$ and $S(\dot H^{-1/2})$.}
$$
\|u\|_{S(\dot H^{1/2})} = \sup_{\substack{(q,r)\; \dot H^{1/2} \text{ admissible} \\
3\leq  r \leq 6^-, \; 4^+ \leq q \leq \infty}} \|u\|_{L_t^qL_x^r} \,
,
$$
where $6^-$ is an arbitrarily preselected and fixed number $<6$;
similarly for $4^+$. We will, in particular, use $(q,r)$ equal to
$(5,5)$, $(20,\frac{10}{3})$, and $(\infty,3)$.  Now we consider
dual Strichartz norms. Let
$$
\|u\|_{S'(L^2)} = \inf_{\substack{(q,r)\; L^2 \text{ admissible} \\
2\leq q \leq \infty, \; 2\leq r \leq 6}} \|u\|_{L_t^{q'}L_x^{r'}},
$$
where $(q',r')$ is the H\"older dual to $(q,r)$.  Also define
$$
\|u\|_{S'(\dot H^{-1/2})} = \inf_{\substack{(q,r)\; \dot H^{-1/2} \text{ admissible} \\
\frac{4}{3}^+ \leq q \leq 2^-, \; 3^+ \leq r \leq 6^-}} \|u\|_{L_t^{q'}L_x^{r'}} \, .
$$

The Strichartz estimates are
$$
\|e^{it\Delta}\phi\|_{S(L^2)} \leq c\|\phi\|_{L^2}
$$
and
$$
\left\| \int_0^t e^{i(t-t')\Delta}f(\cdot,t')dt' \right\|_{S(L^2)} \leq c\|f\|_{S'(L^2)}.
$$
By combining Sobolev embedding with the Strichartz estimates, we obtain
$$\|e^{it\Delta}\phi\|_{S(\dot H^{1/2})} \leq c\|\phi\|_{\dot H^{1/2}}$$
and
\begin{equation}
\label{E:SobStr}
\left\| \int_0^t e^{i(t-t')\Delta}f(\cdot,t') dt'\right\|_{S(\dot H^{1/2})}
\leq c\|D^{1/2}f\|_{S'(L^2)} \,.
\end{equation}
We shall also need the Kato inhomogeneous Strichartz estimate
\cite{Kato} (for further extensions see \cite{F04} and \cite{V})
\begin{equation}
\label{E:Kato}
\left\| \int_0^t e^{i(t-t')\Delta}f(\cdot, t') \, dt' \right\|_{S(\dot H^{1/2})}
\leq c\|f\|_{S'(\dot H^{-1/2})} \, .
\end{equation}
In particular, we will use $L_t^5L_x^5$ and
$L_t^{20}L_x^\frac{10}{3}$ on the left side, and
$L_t^{10/3}L_x^{5/4}$ on the right side.

We extend our notation $S(\dot H^s)$, $S'(\dot H^s)$ as follows: If
a time interval is not specified (that is, if we just write $S(\dot
H^s)$, $S'(\dot H^s)$), then the $t$-norm is evaluated over
$(-\infty,+\infty)$.  To indicate a restriction to a time
subinterval $I\subset (-\infty,+\infty)$, we will write $S(\dot
H^s;I)$ or $S'(\dot H^s; I)$.

\begin{proposition}[Small data]
\label{P:sd} Suppose $\|u_0\|_{\dot H^{1/2}} \leq A$. There is
$\delta_{\textnormal{sd}}=\delta_{\textnormal{sd}}(A)>0$ such that
if $\|e^{it\Delta} u_0\|_{S(\dot H^{1/2})} \leq
\delta_{\textnormal{sd}}$, then $u$ solving \eqref{E:NLS} is global
(in $\dot H^{1/2}$) and
$$
\|u \|_{S(\dot H^{1/2})} \leq 2\,\| e^{it\Delta}u_0\|_{S(\dot
H^{1/2})},
$$
$$
\|D^{1/2}u\|_{S(L^2)} \leq 2\,c\, \|u_0\|_{\dot H^{1/2}}.
$$
(Note that by the Strichartz estimates, the hypotheses are satisfied
if $\|u_0\|_{\dot H^{1/2}} \leq c\delta_{\textnormal{sd}}$.)
\end{proposition}
\begin{proof}
Define
$$
\Phi_{u_0}(v) = e^{it\Delta}u_0
+i\int_0^te^{i(t-t')\Delta}|v|^2v(t')dt'.
$$
Applying the Strichartz estimates, we obtain
$$
\|D^{1/2}\Phi_{u_0}(v)\|_{S(L^2)} \leq c\|u_0\|_{\dot H^{1/2}} + c
\| D^{1/2}(|v|^2v)\|_{L_t^{5/2}L_x^{10/9}}
$$
and
$$
\|\Phi_{u_0}(v)\|_{S(\dot H^{1/2})} \leq \|e^{it\Delta}u_0\|_{S(\dot
H^{1/2})} + c \| D^{1/2}(|v|^2v)\|_{L_t^{5/2}L_x^{10/9}}.
$$
Applying the fractional Leibnitz \cite{KPV93} and H\"older
inequalities
$$
\|D^{1/2}(|v|^2v) \|_{L_t^{5/2}L_x^{10/9}} \leq \|v\|_{L_t^5L_x^5}^2
\|D^{1/2}v\|_{L_t^\infty L_x^2}\leq \|v\|_{S(\dot
H^{1/2})}^2\|D^{1/2}v\|_{S(L^2)}.
$$
Let
$$
\delta_{\text{sd}}\leq   \min\Big( \frac{1}{\sqrt{24}c},
\frac{1}{24cA}\Big).
$$
Then $\Phi_{u_0}: B\to B$, where
$$
B = \left\{ \; v \; \big| \; \|v\|_{S(\dot H^{1/2})} \leq
2\|e^{it\Delta}u_0\|_{S(\dot H^{1/2})}, \; \|D^{1/2}v\|_{S(L^2)}
\leq 2c\|u_0\|_{\dot H^{1/2}} \right\}
$$
and $\Phi_{u_0}$ is a contraction on $B$.
\end{proof}

\begin{proposition}[$H^1$ scattering]
\label{P:persistence} If $u_0\in H^1$,  $u(t)$ is global with
globally finite $\dot H^{1/2}$ Strichartz norm $\|u\|_{S(\dot
H^{1/2})} < +\infty$ and a uniformly bounded $H^1$ norm $\sup_{t\in
[0,+\infty)} \|u(t)\|_{H^1} \leq B$, then $u(t)$ scatters in $H^1$
as $t\to +\infty$.  This means that there exists $\phi^+ \in H^1$
such that
$$
\lim_{t\to +\infty}\|u(t) - e^{it\Delta}\phi^+\|_{H^1} =0.
$$
\end{proposition}
\begin{proof}
Since $u(t)$ solves the integral equation
$$
u(t) = e^{it\Delta}u_0 + i\int_0^t e^{i(t-t')\Delta}(|u|^2u)(t') \,
dt',
$$
we have
\begin{equation}
\label{E:scat} u(t) - e^{it\Delta}\phi^+ = -i\int_t^{+\infty}
e^{i(t-t')\Delta}(|u|^2u)(t')dt',
\end{equation}
where
$$
\phi^+ = u_0 + i\int_0^{+\infty} e^{-it'\Delta}(|u|^2u)(t')dt'.
$$
Applying the Strichartz estimates to \eqref{E:scat}, we have
\begin{align*}
\| u(t) - e^{it\Delta}\phi^+ \|_{H^1} &\leq c\| |u|^2 \,
(1+|\nabla|)u \|_{L_{[t,+\infty)}^{5/2}L_x^{10/9}} \\
&\leq c\|u\|_{L_{[t,+\infty)}^5L_x^5}^2 \|u\|_{L_t^\infty H_x^1}\\
&\leq cB\|u\|_{L_{[t,+\infty)}^5L_x^5}^2 \,.
\end{align*}
Send $t\to +\infty$ in this inequality to obtain the claim.
\end{proof}

The following long-time perturbation result is similar in spirit to
Lemma 3.10 in Colliander-Keel-Staffilani-Takaoka-Tao
\cite{CKSTTAnnals}, although more refined than a direct analogous
version since the smallness condition \eqref{E:longtimesmallness} is
expressed in terms of $S(\dot H^{1/2})$ rather than
$D^{-1/2}S(L^2)$.  This refinement is achieved by employing the Kato
inhomogeneous Strichartz estimates \cite{Kato}.

\begin{proposition}[Long time perturbation theory]
\label{P:longtime} For each $A \gg 1$, there exists
$\epsilon_0=\epsilon_0(A)\ll 1$ and $c=c(A) \gg 1$ such that the following holds.  Let
$u=u(x,t)\in H_x^1$ for all $t$ and solve
$$
i\partial_t u + \Delta u + |u|^2u =0\, .
$$
Let $\tilde u=\tilde u(x,t)\in H_x^1$ for all $t$ and define
$$
e \defeq i\partial_t \tilde u + \Delta \tilde u + |\tilde u|^2 \tilde u \, .
$$
If
$$
\|\tilde u \|_{S(\dot H^{1/2})} \leq  A\, , \quad \| e\|_{S'(\dot
H^{-1/2})} \leq \epsilon_0 \, , \quad \text{and}
$$
\begin{equation}
\label{E:longtimesmallness}
\|e^{i(t-t_0)\Delta} (u(t_0)-\tilde u(t_0)) \|_{S(\dot H^{1/2})} \leq \epsilon_0 \, ,
\end{equation}
then
$$\|u\|_{S(\dot H^{1/2})} \leq c=c(A)<\infty \, .$$
\end{proposition}
\begin{proof}
Let $w$ be defined by $u=\tilde u + w$.  Then $w$ solves the equation
\begin{equation}
\label{E:pert} i\partial_t w + \Delta w + (\tilde u^2\bar w +
2|\tilde u|^2w) + (2\,\tilde u\,|w|^2 + \bar{\tilde u}\,w^2) +|w|^2
w - e=0.
\end{equation}
Since $\|\tilde u\|_{S(\dot H^{1/2})}\leq A$, we can partition
$[t_0,+\infty)$ into $N=N(A)$ intervals\footnote{The number of
intervals depends only on $A$, but the intervals themselves depend
upon the function $\tilde u$.} $I_j=[t_j,t_{j+1}]$ such that for
each $j$, the quantity $\|\tilde u\|_{S(\dot H^{1/2}; I_j)}\leq
\delta$ is suitably small ($\delta$ to be chosen below). The
integral equation version of \eqref{E:pert} with initial time $t_j$ is
\begin{equation}
\label{E:intpert} w(t)= e^{i(t-t_j)\Delta}w(t_j) +i \int_{t_j}^t
e^{i(t-s)\Delta}W(\cdot, s)\, ds,
\end{equation}
where
$$
W = (\tilde u^2\bar w + 2|\tilde u|^2w) + (2\,\tilde u\,|w|^2 +
\bar{\tilde u}\,w^2) +|w|^2w-e.
$$
By applying the Kato Strichartz estimate \eqref{E:Kato} on $I_j$, we
obtain
\begin{equation}
\label{E:55control} \|w\|_{S(\dot H^{1/2}; I_j)}
\begin{aligned}[t]
& \leq \|e^{i(t-t_j)\Delta}w(t_j)\|_{S(\dot H^{1/2};I_j)} + c\,\|\tilde u^2w\|_{L_{I_j}^{10/3}L_x^{5/4}}\\
&+ c \,\| \tilde u w^2 \|_{L_{I_j}^{10/3}L_x^{5/4}}+ c\, \| w^3
\|_{L_{I_j}^{10/3}L_x^{5/4}} +\|e\|_{S'(\dot H^{-1/2}; I_j)} .
\end{aligned}
\end{equation}
Observe
$$
\|\tilde u^2 w\|_{L_{I_j}^{10/3}L_x^{5/4}} \leq \|\tilde
u\|_{L_{I_j}^{20}L_x^{10/3}}^2 \|w\|_{L_{I_j}^5L_x^5}\leq \|\tilde
u\|_{S(\dot H^{1/2}; I_j)}^2 \|w\|_{S(\dot H^{1/2};I_j)} \leq
\delta^2 \|w\|_{S(\dot H^{1/2};I_j)} .
$$
Similarly,
$$
\| \tilde u w^2 \|_{L_{I_j}^{10/3}L_x^{5/4}} \leq \delta \|w\|_{S(\dot H^{1/2};I_j)}^2,
\quad \text{and} \quad \|  w^3 \|_{L_{I_j}^{10/3}L_x^{5/4}}
\leq \|w\|_{S(\dot H^{1/2};I_j)}^3 \, .
$$
Substituting the above estimates in \eqref{E:55control},
\begin{equation}
\label{E:jthbound} \|w\|_{S(\dot H^{1/2}; I_j)}
\begin{aligned}[t]
& \leq \|e^{i(t-t_j)\Delta}w(t_j)\|_{S(\dot H^{1/2};I_j)} + c\delta^2\|w\|_{S(\dot H^{1/2}; I_j)} \\
&+ c\delta \|w\|_{S(\dot H^{1/2}; I_j)}^2 + c\|w\|_{S(\dot H^{1/2};
I_j)}^3+c\epsilon_0 .
\end{aligned}
\end{equation}
Provided
\begin{equation}
\label{E:requirements}
\delta \leq \min\Big( 1, \frac{1}{6c}\Big) \quad \text{and} \quad
\Big(\|e^{i(t-t_j)\Delta}w(t_j)\|_{S(\dot H^{1/2};I_j)}+c\epsilon_0\Big)
\leq \min\Big( 1, \frac{1}{2\sqrt{6c}}\Big) \, ,
\end{equation}
we obtain
\begin{equation}
\label{E:delta} \|w\|_{S(\dot H^{1/2}; I_j)} \leq
2\|e^{i(t-t_j)\Delta}w(t_j)\|_{S(\dot H^{1/2};I_j)} +2c\epsilon_0.
\end{equation}
Now take $t=t_{j+1}$ in \eqref{E:intpert}, and apply $e^{i(t-t_{j+1})\Delta}$
to both sides to obtain
\begin{equation}
\label{E:jthboundprop} e^{i(t-t_{j+1})\Delta}w(t_{j+1})=
e^{i(t-t_j)\Delta}w(t_j) +i \int_{t_j}^{t_{j+1}}
e^{i(t-s)\Delta}W(\cdot, s)\, ds.
\end{equation}
Since the Duhamel integral is confined to $I_j=[t_j,t_{j+1}]$, by
again applying the Kato estimate, similarly to \eqref{E:jthbound} we
obtain the estimate
$$
\|e^{i(t-t_{j+1})\Delta}w(t_{j+1})\|_{S(\dot H^{1/2})}
\begin{aligned}[t]
&\leq \|e^{i(t-t_j)\Delta}w(t_j)\|_{S(\dot H^{1/2})} + c\delta^2\|w\|_{S(\dot H^{1/2}; I_j)} \\
&+ c\delta \|w\|_{S(\dot H^{1/2}; I_j)}^2 + c\|w\|_{S(\dot H^{1/2};
I_j)}^3+c\epsilon_0 .
\end{aligned}
$$
By \eqref{E:delta} and \eqref{E:jthboundprop}, we bound the previous expression to obtain
$$
\|e^{i(t-t_{j+1})\Delta}w(t_{j+1})\|_{S(\dot H^{1/2})} \leq
2\|e^{i(t-t_j)\Delta}w(t_j)\|_{S(\dot H^{1/2})} + 2c\epsilon_0.
$$
Iterating beginning with $j=0$, we obtain
\begin{align*}
\|e^{i(t-t_{j})\Delta}w(t_{j})\|_{S(\dot H^{1/2})}
&\leq 2^j\|e^{i(t-t_0)\Delta}w(t_0)\|_{S(\dot H^{1/2})}+ (2^j-1)2c\epsilon_0 \\
&\leq 2^{j+2}c\epsilon_0 .
\end{align*}
To accommodate the second part of \eqref{E:requirements} for all
intervals $I_j$, $0\leq j \leq N-1$, we require that
\begin{equation}
 \label{E:ep_requirement}
2^{N+2}c\epsilon_0 \leq  \min\Big( 1, \frac{1}{2\sqrt{6c}}\Big) .
\end{equation}
We review the dependence of parameters: $\delta$ is an absolute
constant selected to meet the first part of \eqref{E:requirements}.
We were given $A$, which then determined $N$ (the number of time
subintervals).  The inequality \eqref{E:ep_requirement} specifies
how small $\epsilon_0$ needs to be taken in terms of $N$ (and thus,
in terms of $A$).
\end{proof}

\section{Properties of the ground state}
\label{S:groundstate}

M.\ Weinstein \cite{W83} proved that the sharp constant
$c_{\text{GN}}$ in the Gagliardo-Nirenberg estimate
\begin{equation}
 \label{E:GN}
\|f\|_{L^4}^4 \leq c_{\text{GN}} \|f\|_{L^2}\|\nabla f\|_{L^2}^3
\end{equation}
is attained at the function $Q$ (the ground state described in the
introduction), i.e., $c_{\text{GN}} = \|Q\|_{L^4}^4/(\|Q\|_{L^2}
\|\nabla Q\|_{L^2}^3)$.  By multiplying \eqref{E:Q} by $Q$,
integrating, and applying integration by parts, we obtain
\begin{equation*}
-\|Q\|_{L^2}^2 - \|\nabla Q\|_{L^2}^2 + \|Q\|_{L^4}^4 =0 \, .
\end{equation*}
By multiplying \eqref{E:Q} by $x\cdot \nabla Q$, integrating, and
applying integration by parts, we obtain the Pohozhaev
identity
$$
\frac32\,\|Q\|_{L^2}^2 +\frac12\, \|\nabla Q\|_{L^2}^2 -\frac34\,
\|Q\|_{L^4}^4 =0 \, .
$$
These two identities enable us to obtain the relations
\begin{equation}
\label{E:Qidentities} \|\nabla Q\|_{L^2}^2 = 3\|Q\|_{L^2}^2, \quad
\|Q\|_{L^4}^4 = 4\|Q\|_{L^2}^2 \, ,
\end{equation}
and thus, reexpress
\begin{equation}
 \label{E:GNconstant}
c_{\text{GN}} = \frac{4}{3\|Q\|_{L^2}\|\nabla Q\|_{L^2}} =
\frac{4}{3\sqrt 3 \|Q\|_{L^2}^2}.\footnote{Numerical calculations
show $\|Q\|^2_{L^2(\cR^3)} \cong 18.94$, which gives $c_{\text{GN}} \cong
0.0406$ (in $\cR^3$).}
\end{equation}
We also calculate
\begin{equation}
 \label{E:ME_Q}
M[Q]E[Q] = \|Q\|_{L^2}^2\left(\frac12\|\nabla Q\|_{L^2}^2 -
\frac14\|Q\|_{L^4}^4\right) = \frac16\,\|Q\|_{L^2}^2\|\nabla
Q\|_{L^2}^2= \frac12\,\|Q\|_{L^2}^4 \,.
\end{equation}

For later purposes we recall a version of the Gagliardo-Nirenberg
inequality valid only for radial functions, due to W. Strauss
\cite{S}. In $\mathbb{R}^3$, for any $R>0$, we have
\begin{equation}
 \label{E:Strauss}
\|f\|^4_{L^4(|x|>R)} \leq \frac{c}{R^2} \| f\|^3_{L^2(|x|>R)}
\|\nabla f\|_{L^2(|x|>R)}.
\end{equation}

\section{Global versus blow-up dichotomy}
\label{S:dichotomy}

In this section we show how to obtain Theorem \ref{T:main} part
(\ref{I:global})(a) and part(\ref{I:blow-up}). This was proved in
Holmer-Roudenko \cite{HR06} for general mass supercritical and
energy subcritical NLS equations with $H^1$ initial data, but for
self-containment of this exposition we outline the main ideas here.

Before giving the proof, we observe that the following
quantities are scaling invariant:
$$
\Vert \nabla u \Vert_{L^2} \cdot \Vert u \Vert_{L^2} \quad
\text{and} \quad E[u] \cdot M[u].
$$

Next, we quote a localized version of the virial identity as in
Kenig-Merle \cite{KM06a}.  We refer, for example, to Merle-Rapha\"el
\cite{MR06} or Ozawa-Tsutsumi \cite{OT91} for a proof.
\begin{lemma}[Local virial identity]
 \label{L:localvirial}
Let $\chi \in C_0^{\infty}(\cR^N)$, radially symmetric and $u$ solve
$$
i\partial_t u + \Delta u + |u|^{p-1}u=0.
$$
Then
\begin{equation}
 \label{E:localvirial}
\partial_t^2 \!\!\int \chi(x) \,|u(x,t)|^2 \, dx = 4 \!\!\int \chi'' |\nabla
u|^2 - \int \Delta^2 \chi \, |u|^2 -
4\left(\frac12-\frac1{p+1}\right) \!\!\int \Delta \chi \, |u|^{p+1}.
\end{equation}
\end{lemma}

We prove a slightly stronger version of Theorem \ref{T:main} parts
\eqref{I:global}(a) and \eqref{I:blow-up} that is valid for
nonradial initial condition. The generalization of this theorem to
all mass supercritical and energy critical cases of NLS can be found
in \S \ref{S:extensions} as well as in \cite{HR06}. A different type
of condition for global existence, phrased as $\|u_0\|_{L^2}\leq
\gamma_*(\|\nabla u_0\|_{L^2})$ for a certain monotonic function
$\gamma:\mathbb{R}_+\to \mathbb{R}_+$, is given by B\'egout
\cite{Begout}.

\begin{theorem}[Global versus blow-up dichotomy]
 \label{BlowUp}
Let $u_0 \in {H}^1(\cR^3)$ (possibly non-radial), and let
$I=(-T_*,T^*)$ be the maximal time interval of existence of $u(t)$
solving \eqref{E:NLS}. Suppose that
\begin{equation}
 \label{E:initialME}
M[u_0] \, E[u_0] <  M[Q] \, E[Q].
\end{equation}
If \eqref{E:initialME} holds and
\begin{equation}
 \label{E:less1}
\Vert \nabla u_0 \Vert_{L^2} \Vert u_0 \Vert_{L^2} < \Vert \nabla Q
\Vert_{L^2} \Vert Q \Vert_{L^2}\,,
\end{equation}
then $I = (-\infty, +\infty)$, i.e. the solution exists globally in
time, and for all time $t \in \cR$
\begin{equation}
 \label{E:less2}
\Vert \nabla u(t) \Vert_{L^2} \Vert u_0 \Vert_{L^2} < \Vert \nabla Q
\Vert_{L^2} \Vert Q \Vert_{L^2}.
\end{equation}
If \eqref{E:initialME} holds and
\begin{equation}
 \label{E:greater1}
\Vert \nabla u_0 \Vert_{L^2} \Vert u_0 \Vert_{L^2}
> \Vert \nabla Q \Vert_{L^2} \Vert Q \Vert_{L^2} \,,
\end{equation}
then for $t \in I$
\begin{equation}
 \label{E:greater2}
\Vert \nabla u(t) \Vert_{L^2} \Vert u_0 \Vert_{L^2}
> \Vert \nabla Q \Vert_{L^2}
\Vert Q \Vert_{L^2}.
\end{equation}
Furthermore, if (a) $|x|u_0 \in L^2(\cR^3)$, or (b) $u_0$ is radial,
then $I$ is finite, and thus, the solution blows up in finite time.
\end{theorem}

We recently became aware that the global existence assertion and the
blow-up assertion under the hypothesis  $|x|u_0 \in L^2(\cR^3)$ in
this theorem previously appeared in the literature in
Kuznetsov-Rasmussen-Rypdal-Turitsyn \cite{KRRT}.  We have decided to
keep the proof below since it is short and for the convenience of
the reader (there are significant notational differences between our
paper and theirs).

\begin{remark}\footnote{We thank J. Colliander for supplying this comment.}
Since this theorem applies to the nonradial case, we remark that one
should exploit the Galilean invariance to extend the class of
solutions $u$ to which it applies.  Since $u$ is global
[respectively, blows up in finite time] if and only if a Galilean
transformation of it is global [respectively, blows up in finite
time], given $u$ consider for some $\xi_0\in \mathbb{R}^3$ the
transformed solution
$$
w(x,t) = e^{ix\cdot\xi_0} e^{-it|\xi_0|^2} u(x-2\xi_0t,t).
$$
We compute
$$
\|\nabla w\|_{L^2}^2 = |\xi_0|^2M[u] + 2\xi_0\cdot P[u] + \|\nabla
u\|_{L^2}^2,
$$
where the vector $P[u] = \I \int \bar u \nabla u \, dx$ is the
conserved momentum.  Therefore, $M[w]=M[u]$ and
$$
E[w] = \frac12|\xi_0|^2 M[u] + \xi_0 \cdot P[u] + E[u].
$$
To minimize $E[w]$ and $\|\nabla w\|_{L^2}$, 
we take $\xi_0 = -P[u]/M[u]$. Then we test the
condition \eqref{E:initialME}, and \eqref{E:less1} or
\eqref{E:greater1} for $w$, rather than $u$. This means that for
$P\neq 0$, the hypothesis \eqref{E:initialME} can be sharpened to
$$
M[u]\left( -\frac{P[u]^2}{2M[u]}+E[u]\right) < M[Q]E[Q]
$$
and the hypothesis \eqref{E:less1} can be sharpened to
$$
\left( -\frac{P[u]^2}{M[u]}+ \|\nabla u_0\|_{L^2}^2 \right)
\|u_0\|_{L^2}^2 < \|Q\|_{L^2}^2\|\nabla Q\|_{L^2}^2
$$
and similarly for \eqref{E:greater1}.
\end{remark}

\begin{proof}%
Multiplying the definition of energy by $M[u]$ and using
(\ref{E:GN}), we have
\begin{align*}
M[u] E[u]
&= \frac12 \,\Vert \nabla u \Vert_{L^{2}}^2 \Vert u_0
\Vert^2_{L^2} - \frac1{4} \,\Vert u \Vert^{4}_{L^{4}} \Vert u_0
\Vert^2_{L^2}\\
&\geq \frac12 \,\Vert \nabla u \Vert_{L^{2}}^2 \Vert
u_0 \Vert^2_{L^2}- \frac{1}{4}\,c_{\text{GN}}\, \Vert \nabla u
\Vert_{L^{2}}^{3} \, \Vert u_0 \Vert^3_{L^2}.
\end{align*}
Define $ f(x) = \frac12 \, x^2 - \frac{c_{\text{GN}}}{4} \, x^3$. Then $
f'(x) = x - \frac34 c_{\text{GN}} \, x^2 = x\left (1- \frac34 c_{\text{GN}} \, x
\right)$, and thus, $f'(x) = 0$ when $x_0=0$ and $ x_1 = {\frac43
\frac1{c_{\text{GN}}}} = \|\nabla Q\|_{L^2} \, \|Q\|_{L^2}$ by
\eqref{E:GNconstant}. Note that $f(0)=0$ and $f(x_1) = \frac16 \,
\|\nabla Q\|^2_{L^2} \, \|Q\|^2_{L^2}$. Thus, the graph of $f$ has a
local minimum at $x_0$ and a local maximum at $x_1$. The condition
(\ref{E:initialME}) together with \eqref{E:ME_Q} imply that $M[u_0]
E[u_0] < f(x_1)$. Combining this with energy conservation, we have
\begin{equation}
 \label{E:coercion1}
f(\|\nabla u(t)\|_{L^2} \|u_0\|_{L^2}) \leq M[u_0]\,E[u(t)]=
M[u_0]\, E[u_0]<f(x_1).
\end{equation}

If initially $\|u_0\|_{L^2} \, \Vert \nabla u_0 \Vert_{L^2} < x_1$,
i.e. the condition (\ref{E:less1}) holds, then by
\eqref{E:coercion1} and the continuity of $\|\nabla u(t)\|_{L^2}$ in
$t$, we have $\|u_0\|_{L^2} \,\|\nabla u(t)\|_{L^2} < x_1$ for all
time $t \in I$ which gives (\ref{E:less2}). In particular, the
$\dot{H}^1$ norm of the solution $u$ is bounded, which proves global
existence (and thus, global wellposedness) in this case.

If initially $\|u_0\|_{L^2} \,\Vert \nabla u_0 \Vert_{L^2} > x_1$,
i.e. the condition (\ref{E:greater1}) holds, then by
\eqref{E:coercion1} and the continuity of $\|\nabla u(t)\|_{L^2}$ in
$t$, we have $\|u_0\|_{L^2} \, \|\nabla u(t)\|_{L^2} > x_1$ for all
time $t \in I$ which gives (\ref{E:greater2}). We can refine this
analysis to obtain the following: if the condition
(\ref{E:greater1}) (together with (\ref{E:initialME})) holds, then
there exists $\delta_1>0$ such that $M[u_0]\, E[u_0]< (1-\delta_1)
M[Q]\, E[Q]$, and thus, there exists $\delta_2 = \delta_2(\delta_1)
> 0$ such that $\|u_0\|_{L^2}^2 \, \|\nabla u(t)\|_{L^2}^2
> (1 + \delta_2) \,\|\nabla Q\|^2_{L^2} \, \|Q\|^2_{L^2}$ for all $t
\in I$.

Now if $u$ has a finite variance, we recall the virial identity
$$
\partial_t^2 \int |x|^2 \, |u(x,t)|^2 \, dx = 24 E[u_0] -
4 \Vert \nabla u (t) \Vert^2_{L^2}.
$$
Multiplying both sides by $M[u_0]$ and applying the refinement of
inequalities (\ref{E:initialME}) and (\ref{E:greater2}) mentioned
above as well as (\ref{E:ME_Q}), we get
\begin{align*}
M[u_0] \,\partial_t^2 \!\int \!\!|x|^2 \, |u(x,t)|^2 \, dx & = 24 \,
M[u_0]\, E[u_0] - 4 \Vert \nabla u (t) \Vert^2_{L^2} \,
\Vert u_0 \Vert^2_{L^2}\\
& < 24 \cdot \tfrac16 \, (1-\delta_1) \, \|\nabla Q\|^2_{L^2} \,
\|Q\|^2_{L^2} - 4 (1+\delta_2) \|\nabla Q\|_{L^2}^2 \, \|Q\|_{L^2}^2\\
& = - 4(\delta_1 + \delta_2) \|\nabla Q\|_{L^2}^2 \, \|Q\|_{L^2}^2 < 0,
\end{align*}
and thus, $I$ must be finite, which implies that blow up occurs in
finite time.

If $u_0$ is radial, we use a localized version of the virial
identity \eqref{E:localvirial}.
Choose $\chi(r)$ (radial) such that $\partial_r^2 \chi(r) \leq 2$
for all $r\geq 0$, $\chi(r) = r^2$ for $0\leq r \leq 1$, and
$\chi(r)$ is constant for $r\geq 3$. Let $\chi_m(r) = m^2\chi(r/m)$.
The rest of the argument follows the proof of the main theorem in
Ogawa-Tsutsumi \cite{OT91}, although we include the details here for
the convenience of the reader.  We bound each of the terms in the
local virial identity  \eqref{E:localvirial} as follows, using that
$\Delta\chi_m(r) = 6$ for $r\leq m$ and $\Delta^2\chi_m(r)=0$ for
$r\leq m$:
\begin{align*}
4 \int \chi_m''\,|\nabla u|^2 &\leq 8\int_{\mathbb{R}^3} |\nabla u|^2\,,\\
-\int \Delta^2\chi_m\,|u|^2 &\leq \frac{c}{m^2}\int_{m\leq |x| \leq 3m} |u|^2\,,\\
-\int \Delta \chi_m \,|u|^4 &\leq -6\int_{|x|\leq m} |u|^4 +
c\int_{m\leq |x|\leq 3m} |u|^4 \leq -6\int_{\mathbb{R}^3} |u|^4 +
c'\int_{|x|\geq m} |u|^4\,.
\end{align*}
Adding these three bounds and applying the radial
Gagliardo-Nirenberg estimate \eqref{E:Strauss}, we obtain that for
any large $m > 0$, we have
\begin{align*}
\partial_t^2 \int \chi_m(|x|) \, |u(x,t)|^2 \, dx
&\leq 24 E[u_0] - 4 \int |\nabla u|^2\\
& \qquad + \frac{c_1}{m^{2}} \, \Vert u_0 \Vert_{L^2}^{3} \, \Vert
\nabla u \Vert_{L^2} + \frac{c_2}{m^2} \, \int_{m<|x|} |u|^2.
\end{align*}
Let $\epsilon>0$ be a small constant to be chosen below. Use Young's
inequality in the third term on the right side to separate the
$L^2$-norm and gradient term and then absorb the gradient term into
the second term with the chosen $\ep$. Multiplying the above
expression by $M[u_0]$,
we get
\begin{align}
\notag
\indentalign M[u_0] \, \partial_t^2 \int \chi_m(|x|) \, |u(x,t)|^2 \,dx \\
\notag
&\leq 24 \, E[u_0] M[u_0] - (4 -\ep) \Vert \nabla u \Vert^2_{L^2}
\Vert u_0 \Vert^{2}_{L^2}+ \frac{c(\ep)}{m^{4}} \, \Vert u_0
\Vert_{L^2}^{8} + \frac{c_2}{m^2} \, \Vert u_0 \Vert_{L^2}^{4}\\
\label{E:locvirbd} & \leq - c_3 \|\nabla Q\|_{L^2} \, \|Q\|_{L^2} +
\frac{c(\ep)}{m^{4}} \, \Vert u_0 \Vert_{L^2}^{8} + \frac{c_2}{m^2}
\, \Vert u_0 \Vert_{L^2}^{4},
\end{align}
where
$$
c_3 \equiv  -4 (1-\delta_1)  + (4 -\ep ) (1+\delta_2) = +4\delta_1 -
\epsilon(1+\delta_2).
$$
Select
$\epsilon=\epsilon(\delta_1,\delta_2)>0$ so that $c_3 >0$ and then
take $m=m(c_3, \epsilon, M[u_0])$ large enough so that the right
side of \eqref{E:locvirbd} is bounded by a strictly negative
constant. This implies that the maximal interval of existence $I$ is
finite.
\end{proof}

The next two lemmas provide some additional estimates that hold
under the hypotheses \eqref{E:initialME} and \eqref{E:less1} of
Theorem \ref{BlowUp}.  These estimates will be needed for the
compactness and rigidity results in \S
\ref{S:compactness}-\ref{S:rigidity}.

\begin{lemma}[Lower bound on the convexity of the variance]
Let $u_0 \in H^1(\cR^3)$ satisfy \eqref{E:initialME} and
\eqref{E:less1}. Furthermore, take $\delta> 0$ such that $M[u_0]
E[u_0] < (1-\delta) M[Q]E[Q]$. If $u$ is the solution of the Cauchy
problem \eqref{E:NLS} with initial data $u_0$, then there exists
$c_{\delta} > 0$ such that for all $t \in \cR$
\begin{equation}
 \label{E:comp3}
24E[u]-4\|\nabla u(t)\|_{L^2}^2 =8 \| \nabla u(t) \|_{L^2}^2 - 6
\|u(t)\|_{L^{4}}^{4} \geq c_{\delta} \, \|\nabla u(t) \|_{L^2}^2.
\end{equation}
\end{lemma}
\begin{proof}
By the analysis in the proof of Theorem \ref{BlowUp}, there exists
$\delta_2=\delta_2(\delta)>0$ such that for all $t\in\mathbb{R}$,
\begin{equation}
\label{E:less3} \|u_0\|_{L^2}^2 \|\nabla u(t)\|_{L^2}^2 \leq
(1-\delta_2)^2 \, \|\nabla Q\|_{L^2}^2 \|Q\|_{L^2}^2 .
\end{equation}
Let
$$
h(t) =\frac{1}{\|Q\|_{L^2}^2\|\nabla Q\|_{L^2}^2} \Big(
8\|u_0\|_{L^2}^2\|\nabla u(t)\|_{L^2}^2 -
6\|u_0\|_{L^2}^2\|u(t)\|_{L^4}^4\Big)\,,
$$
and set $g(y) = y^2-y^3$. By the Gagliardo-Nirenberg estimate
(\ref{E:GN}) with sharp constant $c_{\text{GN}}$ and
\eqref{E:GNconstant},
$$
h(t) \geq 8 \, g\Big( \frac{\|\nabla
u(t)\|_{L^2}\|u_0\|_{L^2}}{\|\nabla Q\|_{L^2}\|Q\|_{L^2}} \Big)
\,.
$$
By \eqref{E:less3}, we restrict attention to $0 \leq y \leq
1-\delta_2$.  By an elementary argument, there exists
$c=c(\delta_2)$ such that $g(y)\geq c\,y^2$ if $0\leq y \leq
1-\delta_2$, which completes the proof.
\end{proof}

\begin{lemma}[Comparability of gradient and energy]
\label{L:comparability} Let $u_0 \in H^1(\cR^3)$ satisfy
\eqref{E:initialME} and \eqref{E:less1}. Then
$$
\frac16  \|\nabla u(t) \|_{L^2}^2 \leq E[u] \leq
\frac12  \|\nabla u(t) \|_{L^2}^2 .
$$
\end{lemma}

\begin{proof}
The second inequality is immediate from the definition of energy.
The first one is obtained by observing that
$$
\frac12 \|\nabla u \|^2_{L^2} - \frac14 \|u\|^4_{L^4} \geq \frac12
\|\nabla u\|^2_{L^2} \left(1 - \tfrac12 c_{\text{GN}} \, \|\nabla u\|_{L^2}
\|u\|_{L^2} \right) \geq \frac16 \, \|\nabla u\|^2_{L^2},
$$
where we used \eqref{E:GN}, \eqref{E:GNconstant} and
\eqref{E:less2}.
\end{proof}

In the proofs of Proposition \ref{P:exist_crit} and
\ref{P:crit_compact}, we will need the following result called
\textit{existence of wave operators} since the map $\Omega^+: \psi^+
\mapsto v_0$ is called the \textit{wave operator} (see the
proposition for the meaning of $\psi^+$ and $v_0$).

\begin{proposition}[Existence of wave operators]
\label{P:waveop}
Suppose $\psi^+\in H^1$ and
\begin{equation}
\label{E:psi_plus} \frac12\|\psi^+\|_{L^2}^2 \|\nabla
\psi^+\|_{L^2}^2 < M[Q]E[Q] .
\end{equation}
Then there exists $v_0\in H^1$ such that $v$ solving \eqref{E:NLS}
with initial data $v_0$ is global in $H^1$ with
$$
\|\nabla v(t)\|_{L^2}\|v_0\|_{L^2} \leq \|Q\|_{L^2}\|\nabla
Q\|_{L^2}, \quad M[v] = \|\psi^+\|_{L^2}^2, \quad E[v] =
\frac12\|\nabla \psi^+\|_{L^2}^2,
$$
and
$$
\lim_{t\to +\infty} \|v(t)-e^{it\Delta}\psi^+\|_{H^1}=0 .
$$
Moreover, if $\|e^{it\Delta}\psi^+\|_{S(\dot H^{1/2})} \leq
\delta_{\textnormal{sd}}$, then
$$
\|v_0\|_{\dot H^{1/2}} \leq 2\,\|\psi^+\|_{\dot H^{1/2}} \quad
\text{and} \quad \| v\|_{S(\dot H^{1/2})} \leq
2\,\|e^{it\Delta}\psi^+\|_{S(\dot H^{1/2})}.
$$
\end{proposition}

\begin{proof}
We want to solve the integral equation
\begin{equation}
\label{E:v_int_eq}
v(t) = e^{it\Delta}\psi^+ - i\int_t^{+\infty} e^{i(t-t')\Delta}(|v|^2v)(t') \, dt'
\end{equation}
first for $t\geq T$ with $T$ large. This is achieved as in the proof
of the small data scattering theory (Proposition \ref{P:sd}), since
there exists $T \gg 0$ such that $\|e^{it\Delta}\psi^+\|_{S(\dot
H^{1/2}; [T,+\infty))} \leq \delta_{\text{sd}}$. By estimating
\eqref{E:v_int_eq}, we obtain
\begin{align*}
\|\nabla v\|_{S(L^2; [T,+\infty))} &\leq c\|\psi^+\|_{\dot H^1} + c\|\nabla (v^3)\|_{S'(L^2)} \\
&\leq c\|\psi^+\|_{\dot H^1} + c\|\nabla
v\|_{S(L^2;[T,+\infty))}\|v\|_{S(\dot H^{1/2};[T,+\infty))}^2 ,
\end{align*}
where in the last step, we used $\|\cdot \|_{S'(L^2)}\leq
\|\cdot\|_{L_t^{10/7}L_x^{10/7}}$ and the H\"older partition
$\frac{7}{10}=\frac{3}{10}+\frac{1}{5}+\frac{1}{5}$.  Thus,
$\|\nabla v\|_{S(L^2; [T,+\infty))} \leq 2\,c\,\|\psi^+\|_{\dot
H^1}$. Using this, we obtain similarly,
$$
\|\nabla (v(t)-e^{it\Delta}\psi^+)\|_{S(L^2;[T,+\infty))} \to 0 \;
\text{as} \; T\to +\infty .
$$
Since $v(t)-e^{it\Delta}\psi^+ \to 0$ in $H^1$ as $t\to +\infty$,
$e^{it\Delta}\psi^+\to 0$ in $L^4$ as $t\to +\infty$, and $\|\nabla
e^{it\Delta}\psi^+\|_{L^2}$ is conserved, we have
$$
E[v] = \lim_{t\to+\infty} \Big(\frac12\|\nabla e^{it\Delta}\psi^+
\|_{L^2}^2- \frac14\|e^{it\Delta}\psi^+\|_{L^4}^4 \Big) = \frac12\|\psi^+\|_{L^2}^2 \,.
$$
Immediately, we obtain $M[v] = \|\psi^+\|_{L^2}^2$. Note that we now
have $M[v]E[v] < M[Q]E[Q]$ by \eqref{E:psi_plus}. Observe that
\begin{align*}
\lim_{t\to +\infty} \|\nabla v(t)\|_{L^2}^2 \|v\|_{L^2}^2
&= \lim_{t \to +\infty} \|\nabla e^{it\Delta}\psi^+\|_{L^2}^2 \| e^{it\Delta}\psi^+\|_{L^2}^2 \\
&= \|\nabla \psi^+\|_{L^2}^2 \| \psi^+\|_{L^2}^2\\
&\leq 2M[Q]E[Q]\\
&= \frac13\|\nabla Q\|_{L^2}^2 \|Q\|_{L^2}^2,
\end{align*}
where in the last two steps we used \eqref{E:psi_plus} and
\eqref{E:ME_Q}. Take $T$ sufficiently large so that $\|\nabla
v(T)\|_{L^2}\|v\|_{L^2} \leq \|\nabla Q\|_{L^2} \|Q\|_{L^2}$.  By
Theorem \ref{BlowUp}, we can evolve $v(t)$ from
$T$ back to time $0$.
\end{proof}

\section{Compactness}
\label{S:compactness}

\begin{definition}
Suppose $u_0\in H^1$ and $u$ is the corresponding $H^1$ solution to
\eqref{E:NLS} and $[0,T^*)$ the maximal forward time interval of
existence. We shall say that $\textnormal{SC}(u_0)$ holds if
$T^*=+\infty$ and $\|u\|_{S(\dot H^{1/2})} <\infty$.
\end{definition}

To prove Theorem \ref{T:main}\eqref{I:global}(b), we must show that
if $\|u_0\|_{L^2} \|\nabla u_0\|_{L^2} < \|Q\|_{L^2}\|\nabla
Q\|_{L^2}$, and $M[u]E[u] < M[Q]E[Q]$, then $\text{SC}(u_0)$ holds.
We already know that under these hypotheses, we have an \textit{a
priori} bound on $\|\nabla u(t)\|_{L^2}$, and thus, the maximal
forward time of existence is $T^*=+\infty$ (this is the energy
\textit{subcritical} case). The goal is, therefore, to show that the
global-in-time $\dot H^{1/2}$ Strichartz norm is finite.

By Lemma \ref{L:comparability}, if $M[u]E[u] <
\frac16\,\delta_{\text{sd}}^4$ and $\|u_0\|_{L^2} \|\nabla u_0\|_{L^2}
< \|Q\|_{L^2}\|\nabla Q\|_{L^2}$, then
$$
\|u_0\|_{\dot H^{1/2}}^4 \leq \|u_0\|_{L^2}^2\|\nabla u_0\|_{L^2}^2 \leq
6M[u]E[u] \leq \delta_{\text{sd}}^4 .
$$
Thus, by the small data theory (Proposition \ref{P:sd}),
$\text{SC}(u_0)$ holds. Let $(ME)_\text{c}$ be the number defined as
the supremum over all $\delta$ for which the following statement
holds true:  ``If $u_0$ is radial with $\|u_0\|_{L^2} \|\nabla
u_0\|_{L^2} < \|Q\|_{L^2}\|\nabla Q\|_{L^2}$ and $M[u]E[u] <
\delta$, then $\text{SC}(u_0)$ holds.'' We then clearly have
$0<\frac16\,\delta_{\text{sd}}^4 \leq (ME)_\text{c}$. If
$(ME)_\text{c} \geq M[Q]E[Q]$, then Theorem
\ref{T:main}\eqref{I:global}(b) is true.  We, therefore, proceed
with the proof of Theorem \ref{T:main}\eqref{I:global}(b) by
assuming that $(ME)_\text{c} < M[Q]E[Q]$ and ultimately deduce a
contradiction (much later, in \S\ref{S:rigidity}). By definition of
$(ME)_\text{c}$, we have
\begin{enumerate}
\item[(C.1)]
If $u_0$ is radial and $\|u_0\|_{L^2}\|\nabla u_0\|_{L^2} <
\|Q\|_{L^2}\|\nabla Q\|_{L^2}$ and $M[u]E[u]< (ME)_\text{c}$, then
$\text{SC}(u_0)$ holds.
\item[(C.2)]
There exists a sequence of radial solutions $u_n$ to \eqref{E:NLS}
with corresponding initial data $u_{n,0}$ such that
$\|u_{n,0}\|_{L^2}\|\nabla u_{n,0}\|_{L^2}< \|Q\|_{L^2}\|\nabla
Q\|_{L^2}$ and $M[u_n]E[u_n] \searrow (ME)_\text{c}$ as $n\to
+\infty$, for which $\text{SC}(u_{n,0})$ does not hold for any $n$.
\end{enumerate}

\newcommand{\crit}{{\textnormal{c}}}
\newcommand{\crin}{{\textnormal{c},0}}

The goal of this section is to use the above sequence $u_{n,0}$
(rescaled so that $\|u_{n,0}\|_{L^2}=1$ for all $n$) to prove the
existence of an $H^1$ radial solution $u_\text{c}$ to \eqref{E:NLS}
with initial data $u_{\text c,0}$ such that $\|u_\crin\|_{L^2}
\|\nabla u_\crin\|_{L^2}< \|Q\|_{L^2}\|\nabla Q\|_{L^2}$ and
$M[u_\crit ]E[u_\crit ]= (ME)_\text{c}$ for which
$\text{SC}(u_\crin)$ does not hold (Proposition \ref{P:exist_crit}).
Moreover, we will show that $K = \{ \, u_\crit(t) \, | \, 0\leq t <
+\infty \, \}$ is precompact in $H^1$ (Proposition
\ref{P:crit_compact}), which will enable us to show that for each
$\epsilon>0$, there is an $R> 0$ such that, uniformly in $t$, we
have
$$
\int_{|x|>R} |\nabla u_\crit(t,x)|^2 dx \leq \epsilon
$$
(Lemma \ref{L:unif_small}). This will then play into the
rigidity theorem of the next section that will ultimately lead to a
contradiction.

Before stating and proving Proposition \ref{P:exist_crit}, we
introduce some preliminaries in the spirit of the results of Keraani
\cite{K01}, since we are not able to directly apply his results as
was possible in Kenig-Merle \cite{KM06a}.  Note in the following
lemma that $\phi_n$, $\psi^j$ and $W^M_n$ are functions of $x$
alone, in notational contrast to the analogous lemma in Keraani
(\cite{K01} Proposition 2.6).

\begin{lemma}[Profile expansion]
 \label{L:pe}
Let $\phi_n(x)$ be a radial uniformly bounded sequence in $H^1$.
Then for each $M$ there exists a subsequence of $\phi_n$, also
denoted $\phi_n$, and
\begin{enumerate}
\item
For each $1\leq j \leq M$, there exists a (fixed in $n$) radial
profile $\psi^j(x)$ in $H^1$,
\item
For each $1\leq j\leq M$, there exists a sequence (in $n$) of time shifts $t^j_n$,
\item
There exists a sequence (in $n$) of remainders $W_n^M(x)$ in $H^1$,
\end{enumerate}
such that
$$
\phi_n = \sum_{j=1}^M e^{-it^j_n\Delta}\psi^j + W^M_n.
$$
The time sequences have a pairwise divergence property: For $1\leq i
\neq j \leq M$, we have
\begin{equation}
 \label{E:pwdiverge}
\lim_{n\to +\infty} |t^i_n - t^j_n| = +\infty .
\end{equation}
The remainder sequence has the following asymptotic smallness
property\footnote{We can always pass to a subsequence in $n$ with the
property that $\| e^{it\Delta}W^M_n\|_{S(\dot H^{1/2})}$ converges.
Therefore, we use $\lim$ and not $\limsup$ or $\liminf$. Similar
remarks apply for the limits that appear in the Pythagorean
expansion.}:
\begin{equation}
 \label{E:remainder_small}
\lim_{M\to +\infty} \left[ \lim_{n\to +\infty} \|
e^{it\Delta}W^M_n\|_{S(\dot H^{1/2})} \right] = 0 .
\end{equation}
For fixed $M$ and any $0\leq s \leq 1$, we have the asymptotic Pythagorean expansion
\begin{equation}
\label{E:pythag} \|\phi_n\|_{\dot H^s}^2 = \sum_{j=1}^M
\|\psi^j\|_{\dot H^s}^2 + \|W^M_n\|_{\dot H^s}^2 + o_n(1).
\end{equation}
\end{lemma}

Note that we do \textit{not} claim that the remainder $W^M_n$ is
small in any Sobolev norm, i.e.\ for all we know it \textit{might}
be true that for some $s$, $0\leq s\leq 1$, we have
$$
\liminf_{M\to +\infty} \left[ \lim_{n\to +\infty} \|W^M_n\|_{\dot H^s} \right] > 0 \, .
$$
Fortunately, the Strichartz norm smallness \eqref{E:remainder_small}
will suffice in our application.

\begin{proof}
Since $\phi_n$ is assumed uniformly bounded in $H^1$, let $c_1$ be
such that $\|\phi_n\|_{H^1} \leq c_1$.

Note the interpolation inequality
$$
\|v \|_{L_t^qL_x^r} \leq \|v\|_{L_t^{\tilde q} L_x^{\tilde
r}}^{1-\theta} \| v \|_{L_t^\infty L_x^3}^\theta \, ,
$$
where $(q,r)$ is any $\dot H^{1/2}$ Strichartz admissible pair (so
$\frac2q+\frac3r=1$ and $3\leq r <+\infty$), $\theta =
\frac{3}{2r-3}$ (so $0<\theta \leq 1$), $\tilde r =2r$, and $\tilde
q = \frac{4r}{2r-3}$.  Observe that $(\tilde q, \tilde r)$ is also
$\dot H^{1/2}$ Strichartz admissible.  By this inequality and the
Strichartz estimates (for $0<\theta< \tfrac35$)\footnote{This
restriction is for consistency with our definition of $S(\dot
H^{1/2})$ in \S \ref{S:localtheory}.}, we get
$$
\|e^{it\Delta}W_n^M \|_{L_t^q L_x^r} \leq c\,\|W_n^M\|_{\dot
H^{1/2}}^{1-\theta}\|e^{it\Delta} W_n^M\|_{L_t^\infty L_x^3}^\theta.
$$
Since we will have $\|W_n^M\|_{\dot H^{1/2}} \leq c_1$, it will
suffice for us to show that
$$
\lim_{M\to +\infty} \left[ \limsup_{n\to +\infty}\|e^{it\Delta}
W_n^M\|_{L_t^\infty L_x^3}\right] = 0 \, .
$$

Let $A_1 \equiv \limsup_{n\to +\infty}
\|e^{it\Delta}\phi_n\|_{L_t^\infty L_x^3}$.  If $A_1=0$, the proof
is complete with $\psi^j=0$ for all $1 \leq j \leq M$.  Suppose
$A_1>0$.  Pass to a subsequence so that $\lim_{n\to +\infty}
\|e^{it\Delta}\phi_n\|_{L_t^\infty L_x^3}= A_1$.  We will show that
there is a time sequence $t_n^1$ and a profile $\psi^1\in H^1$ such
that $e^{it^1_n\Delta} \phi_n \rightharpoonup \psi^1$ and
$\|\psi^1\|_{\dot H^{1/2}} \geq \frac{A_1^5}{2^{10}c_1^4}$.    For
$r>1$ yet to be chosen, let $\chi(x)$ be a radial Schwartz function
such that $\hat \chi(\xi) = 1$ for $\frac1r \leq |\xi| \leq r$ and
$\hat \chi(\xi)$ is supported in  $\frac1{2r} \leq |\xi| \leq 2r$.

By Sobolev embedding,
\begin{align*}
\|e^{it\Delta}\phi_n - \chi * e^{it\Delta} \phi_n\|^2_{L_t^\infty L_x^3} %
&\leq \int |\xi| (1-\hat \chi(\xi))^2 |\hat \phi_n(\xi)|^2 \, d\xi \\
&\leq \int_{|\xi|\leq \frac1r} |\xi| |\hat \phi_n(\xi)|^2 \, d\xi +  \int_{|\xi|\geq r} |\xi| |\hat \phi_n(\xi)|^2 \, d\xi\\
& \leq \frac1r \|\phi_n\|_{L^2}^2 + \frac1r\|\phi_n\|_{\dot H^1}^2\\
&\leq \frac{c_1^2}{r} \, .
\end{align*}
Take $r=\frac{16c_1^2}{A_1^2}$ so that $\frac{c_1^2}{r} =
\frac{A_1^2}{16}$, and then we have for $n$ large
$$
\|\chi * e^{it\Delta}\phi_n\|_{L_t^\infty L_x^3} \geq \frac12A_1 \,.
$$
Note that
\begin{align*}
\|\chi * e^{it\Delta}\phi_n\|_{L_t^\infty L_x^3}^3 &\leq \|\chi * e^{it\Delta}\phi_n\|_{L_t^\infty L_x^2}^2\|\chi * e^{it\Delta}\phi_n\|_{L_t^\infty L_x^\infty} \\
&\leq \|\phi_n\|_{L^2}^2\|\chi * e^{it\Delta}\phi_n\|_{L_t^\infty
L_x^\infty},
\end{align*}
and thus, we have
$$
\|\chi * e^{it\Delta}\phi_n\|_{L_t^\infty L_x^\infty} \geq
\frac{A_1^3}{8c_1^{2}} .
$$
Since $\phi_n$ are radial functions, so are $\chi *
e^{it\Delta}\phi_n$, and by the radial Gagliardo-Nirenberg
inequality, we obtain
$$
\| \chi * e^{it\Delta} \phi_n \|_{L_t^\infty L^\infty_{\{|x|\geq
R\}}} \leq \frac{1}{R} \, \|\chi * e^{it\Delta}
\phi_n\|_{L_x^2}^{1/2} \, \|\nabla \chi * e^{it\Delta}
\phi_n\|_{L_x^2}^{1/2} \leq \frac{c_1}{R} .
$$
Therefore, by selecting $R$ large enough
$$
\|\chi * e^{it\Delta}\phi_n\|_{L_t^\infty L^\infty_{\{|x|\leq R\}}}
\geq \frac{A_1^3}{16c_1^{2}} .
$$
Let $t^1_n$ and $x^1_n$ (with $|x^1_n|\leq R$) be sequences such
that for each $n$,
$$
|\chi * e^{it^1_n\Delta}\phi_n (x^1_n)| \geq
\frac{A_1^3}{32c_1^{2}},
$$
or, written out,
$$
\left|\int_{\mathbb{R}^3} \chi(x_n^1-y) \,e^{it^1_n\Delta}\phi_n
(y)\, dy \right| \geq \frac{A_1^3}{32c_1^{2}} .
$$
Pass to a subsequence such that $x_n^1 \to x^1$ (possible since
$|x^1_n|\leq R$). Then since $\|\chi(x^1-\cdot) -
\chi(x_n^1-\cdot)\|_{L^{3/2}} \to 0$ as $n\to +\infty$, we have
$$
\left|\int_{\mathbb{R}^3} \chi(x^1-y) \,e^{it^1_n\Delta}\phi_n (y)\,
dy \right| \geq \frac{A_1^3}{64c_1^{2}} .
$$
Consider the sequence $e^{it^1_n\Delta}\phi_n$, which is uniformly
bounded in $H^1$. Pass to a subsequence so that
$e^{it^1_n\Delta}\phi_n \rightharpoonup \psi^1$, with $\psi^1\in
H^1$ radial and $\|\psi^1\|_{H^1} \leq \limsup \|\phi_n\|_{H^1} \leq
c_1$.  By the above, we have
$$
\left|\int_{\mathbb{R}^3} \chi(x^1-y) \,\psi^1 (y)\, dy \right| \geq
\frac{A_1^3}{64c_1^{2}} .
$$
By Plancherel and Cauchy-Schwarz applied to the left side of the
above inequality, we obtain
$$
\|\chi\|_{\dot H^{-1/2}} \|\psi^1\|_{\dot H^{1/2}} \geq
\frac{A_1^3}{64c_1^{2}} .
$$
By converting to radial coordinates, we can estimate $\|\chi\|_{\dot
H^{-1/2}}  \leq r$.  Thus,
$$
\|\psi^1\|_{\dot H^{1/2}} \geq \frac{A_1^3}{64c_1^{2}}\cdot
\frac{1}{r} = \frac{A_1^5}{2^{10}c_1^4}.
$$
Let $W^1_n = \phi_n - e^{-it^1_n\Delta}\psi^1$.  Since
$e^{it^1_n\Delta}\phi_n \rightharpoonup \psi^1$, we have that for
any $0\leq s \leq 1$
\begin{equation}
 \label{E:psi1}
\la \phi_n, e^{-it^1_n\Delta}\psi^1 \ra_{\dot H^s} = \la
e^{it_n^1\Delta}\phi_n, \psi^1 \ra_{\dot H^s} \to \|\psi^1\|_{\dot
H^s}^2 ,
\end{equation}
and, by expanding $\|W^1_n\|_{\dot H^s}^2$, we obtain
$$
\lim_{n\to +\infty} \|W^1_n\|_{\dot H^s}^2 =
\lim_{n\to+\infty}\|\phi_n\|_{\dot H^s}^2 - \|\psi^1\|_{\dot H^s}^2.
$$
From this with $s=1$ and $s=0$ we deduce that $\|W^1_n\|_{H^1}
\leq c_1$.

Let $A_2 = \limsup_{n\to +\infty} \|e^{it\Delta}W^1_n\|_{L_t^\infty
L_x^3}$. If $A_2 =0$, then we are done.  If $A_2>0$, then repeat the
above argument, with $\phi_n$ replaced by $W^1_n$ to obtain a
sequence of time shifts $t^2_n$ and a profile $\psi^2\in H^1$ such
that $e^{it^2_n\Delta}W^1_n\rightharpoonup \psi^2$ and
$$
\|\psi^2\|_{\dot H^{1/2}} \geq \frac{A_2^5}{2^{10}c_1^4}.
$$
We claim that $|t^2_n-t^1_n| \to +\infty$.  Indeed, suppose we pass
to a subsequence such that $t^2_n-t^1_n \to t^{21}$ finite.  Then
$$
e^{i(t^2_n-t^1_n)\Delta}[e^{it_n^1\Delta}\phi_n - \psi^1] =
e^{it^2_n\Delta}[\phi_n - e^{-it^1_n\Delta}\psi^1]
=e^{it^2_n\Delta}W^1_n \rightharpoonup \psi^2.
$$
Since $t^2_n-t^1_n \to t^{21}$ and $e^{it_n^1\Delta}\phi_n - \psi^1
\rightharpoonup 0$, the left side of the above expression converges
weakly to $0$, so $\psi^2=0$, a contradiction.  Let $W^2_n = \phi_n
- e^{it^1_n\Delta}\psi^1 - e^{it^2_n\Delta}\psi^2$.  Note that
\begin{align*}
\la \phi_n, e^{it^2_n\Delta}\psi^2 \ra_{\dot H^s}
&= \la e^{-it_n^2\Delta}\phi_n, \psi^2 \ra_{\dot H^s} \\
&= \la e^{-it_n^2\Delta}(\phi_n- e^{it^1_n\Delta}\psi^1), \psi^2 \ra_{\dot H^s} + o_n(1)\\
&=\la e^{-it_n^2\Delta}W^1_n, \psi^2 \ra_{\dot H^s} + o_n(1) \\
&\to \|\psi^2\|_{\dot H^s}^2,
\end{align*}
where the second line follows from the fact that $|t^1_n-t^2_n|\to
\infty$.  \ Using this and \eqref{E:psi1}, we compute
$$
\lim_{n\to+\infty} \|W^2_n\|_{\dot H^s}^2 = \lim_{n\to +\infty}
\|\phi_n\|_{\dot H^s}^2 - \|\psi^1\|_{\dot H^s}^2 - \|\psi^2\|_{\dot
H^s}^2,
$$
and thus, $\|W^2_n\|_{H^1} \leq c_1$.

We continue inductively, constructing a sequence $t^M_n$ and a
profile $\psi^M$ such that $e^{it^M_n\Delta}W^{M-1}_n\rightharpoonup
\psi^M$ and
\begin{equation}
\label{E:psibd} \|\psi^M\|_{\dot H^{1/2}} \geq
\frac{A_M^5}{2^{10}c_1^4}.
\end{equation}
Suppose $1\leq j < M$. We shall show that $|t^M_n - t^j_n| \to
+\infty$ inductively by assuming that $|t^M_n-t^{j+1}_n| \to
+\infty, \ldots, |t^M_n-t^{M-1}_n| \to +\infty$. Suppose, passing to
a subsequence that $t^M_n - t^j_n \to t^{Mj}$ finite. Note that
$$
e^{i(t_n^M-t^j_n)\Delta}(e^{it^j_n\Delta}W^{j-1}_n -\psi^j) -
e^{i(t^M_n-t^{j+1}_n)\Delta}\psi^{j+1} - \cdots -
e^{i(t_n^M-t^{M-1}_n)\Delta}\psi^{M-1} = e^{it^M_n\Delta}
W^{M-1}_n.
$$
The left side converges weakly to $0$, while the right
side converges weakly to $\psi^M$, which is nonzero; contradiction.
This proves \eqref{E:pwdiverge}. Let $W^M_n = \phi_n -
e^{-it^1_n\Delta}\psi^1 - \cdots - e^{-it^M_n\Delta}\psi^M$.  Note
that
\begin{align*}
\la \phi_n, e^{-it_n^M\Delta}\psi^M \ra &= \la e^{it^M_n\Delta}\phi_n, \psi^M \ra_{\dot H^s}\\
&= \la e^{it^M_n\Delta}(\phi_n- e^{it^1_n\Delta}\psi^1 - \cdots - e^{it^{M-1}_n\Delta}\psi^{M-1}), \psi^M \ra_{\dot H^s} + o_n(1)\\
&= \la e^{it^M_n\Delta}W^{M-1}_n, \psi^M \ra_{\dot H^s} + o_n(1),
\end{align*}
where the middle line follows from the pairwise divergence property
\eqref{E:pwdiverge}. Thus, $\la \phi_n, e^{-it_n^M\Delta}\psi^M \ra
\to \|\psi^M\|_{\dot H^s}^2$. The expansion \eqref{E:pythag} is then
shown to hold by expanding $\|W^M_n\|_{\dot H^s}^2$.

By \eqref{E:psibd} and \eqref{E:pythag} with $s=\frac12$, we have
$$
\sum_{M=1}^{+\infty}  \left(\frac{A_M^5}{2^{10}c_1^4}\right)^2  \leq
\lim_{n\to +\infty} \|\phi_n\|_{\dot H^{1/2}}^2 \leq c_1^2,
$$
and hence, $A_M\to 0$ as $M\to +\infty$.
\end{proof}

\begin{corollary}[Energy Pythagorean expansion]
\label{C:energy_expand}
In the situation of Lemma \ref{L:pe}, we have
\begin{equation}
\label{E:en_expand} E[\phi_n] = \sum_{j=1}^M
E[e^{-it^j_n\Delta}\psi^j] + E[W^M_n] + o_n(1) .
\end{equation}
\end{corollary}
\begin{proof}
We will use the compact embedding $H^1_{\text{rad}} \hookrightarrow
L^4_{\text{rad}}$ (which follows from the radial Gagliardo-Nirenberg
estimate of Strauss \cite{S}) to address a $j$ for which $t^j_n$
converges to a finite number (if one exists). We will also use the
decay of linear Schr\"odinger solutions in the $L^4$ norm as time
$\to \infty$.

There are two cases to consider.

\textit{Case 1}. There exists some $j$ for which $t^j_n$ converges
to a finite number, which without loss we assume is $0$.  In this
case we will show that
$$
\lim_{n\to +\infty} \|W_n^M\|_{L_x^4} = 0, \quad \text{for }M>j,
$$
$$
\lim_{n\to +\infty}  \|e^{-it_n^i\Delta}\psi^i\|_{L_x^4}=0, \quad
\text{for all } i \neq j,
$$
and
$$
\lim_{n\to +\infty} \|\phi_n\|_{L^4} = \|\psi^j\|_{L^4},
$$
which, combined with \eqref{E:pythag} for $s=1$, gives
\eqref{E:en_expand}.

\textit{Case 2}.  For all $j$, $|t^j_n| \to \infty$. In this case we
will show that
$$
\lim_{n\to +\infty}  \|e^{-it_n^j\Delta}\psi^j\|_{L_x^4}=0, \quad \text{for all }j
$$
and
$$
\lim_{n\to +\infty} \|\phi_n\|_{L^4} =
\lim_{n\to+\infty}\|W^M_n\|_{L^4},
$$
which, combined with \eqref{E:pythag} for $s=1$, gives \eqref{E:en_expand}.

\textit{Proof of Case 1}. In this situation, we have, from the proof
of Lemma \ref{L:pe} that  $W_n^{j-1}\rightharpoonup \psi^j$.  By the
compactness of the embedding $H_{\text{rad}}^1 \hookrightarrow
L_{\text{rad}}^4$, it follows that $W_n^{j-1}\to \psi^j$ strongly in
$L^4$.  Let $i\neq j$.  Then we claim that
$\|e^{it_n^i\Delta}\psi^i\|_{L^4} \to 0$ as $n\to \infty$.  Indeed,
since $t_n^j=0$, by \eqref{E:pwdiverge}, we have $|t_n^i|\to
+\infty$.  For a function $\tilde \psi^i\in \dot H^{3/4}\cap
L^{4/3}$, from Sobolev embedding and the $L^p$ spacetime decay
estimate of the linear flow, we obtain
$$
\|e^{it_n^j\Delta}\psi^i\|_{L^4} \leq c\|\psi^i-\tilde
\psi^i\|_{\dot H^{3/4}} + \frac{c}{|t_n^i|^{1/4}}\|\tilde
\psi^i\|_{L^{4/3}}.
$$
By approximating $\psi^i$ by $\tilde\psi^i\in C_c^\infty$ in $\dot
H^{3/4}$ and sending $n\to +\infty$, we obtain the claim. Recalling
that
$$
W_n^{j-1}=\phi_n - e^{-it_n^1\Delta}\psi^1 - \cdots -
e^{-it_n^{j-1}\Delta}\psi^{j-1},
$$
we conclude that $\phi_n\to \psi^j$ strongly in $L^4$. Recalling
that
$$
W_n^M=(W_n^{j-1}-\psi^j) - e^{-it_n^{j+1}\Delta}\psi^{j+1} - \cdots
- e^{-it_n^M\Delta}\psi^M,
$$
we also conclude that $W_n^M \to 0$ strongly in $L^4$ for $M>j$.

\textit{Proof of Case 2}.  Similar to the proof of Case 1.
\end{proof}

\begin{proposition}[Existence of a critical solution]
\label{P:exist_crit} There exists a global $(T^*=+\infty)$ solution
$u_\crit$ in $H^1$ with initial data $u_\crin$ such that
$\|u_\crin\|_{L^2}=1$,
$$
E[u_\crit] = (ME)_\crit < M[Q]E[Q],
$$
$$
\|\nabla u_\crit(t)\|_{L^2} < \|Q\|_{L^2}\|\nabla Q\|_{L^2} \quad
\text{for all } 0\leq t<+\infty,
$$
and
$$
\|u_\crit \|_{S(\dot H^{1/2})} = +\infty.
$$
\end{proposition}
\begin{proof}
We consider the sequence $u_{n,0}$ described in the introduction to
this section. Rescale it so that $\|u_{n,0}\|_{L^2} =1$; this
rescaling does not affect the quantities $M[u_n]E[u_n]$ and
$\|u_{n,0}\|_{L^2}\|\nabla u_{n,0}\|_{L^2}$.  After this rescaling,
we have $\|\nabla u_{n,0}\|_{L^2} < \|Q\|_{L^2}\|\nabla Q\|_{L^2}$
and $E[u_n] \searrow (ME)_\text{c}$.  Each $u_n$ is global and
non-scattering, i.e.\ $\|u_n\|_{S(\dot H^{1/2})} = +\infty$.  Apply
the profile expansion lemma (Lemma \ref{L:pe}) to $u_{n,0}$ (which
is now uniformly bounded in $H^1$) to obtain
\begin{equation}
\label{E:profexpand} u_{n,0} = \sum_{j=1}^M e^{-it_n^j\Delta}\psi^j
+ W_n^M ,
\end{equation}
where $M$ will be taken large later. By the energy Pythagorean
expansion (Corollary \ref{C:energy_expand}), we have
$$
\sum_{j=1}^M \lim_{n\to+\infty} E[e^{-it_n^j\Delta}\psi^j] +
\lim_{n\to +\infty} E[W_n^M] = \lim_{n\to +\infty}E[u_{n,0}]=
(ME)_\text{c},
$$
and thus (recalling that each energy is $\geq 0$ -- see Lemma
\ref{L:comparability}),
\begin{equation}
\label{E:energy-below} \lim_{n\to+\infty} E[e^{-it_n^j\Delta}\psi^j]
\leq (ME)_\text{c} \quad \forall \; j .
\end{equation}
Also by $s=0$ of \eqref{E:pythag}, we have
\begin{equation}
\label{E:mass-sum} \sum_{j=1}^M M[\psi^j] +
\lim_{n\to+\infty}M[W_n^M] = \lim_{n\to+\infty} M[u_{n,0}] = 1.
\end{equation}
Now we consider two cases; we will show that Case 1 leads to a
contradiction and thus does not occur; Case 2 will manufacture the
desired critical solution $u_\text{c}$.

\textit{Case 1}.  More than one $\psi^j \neq 0$. By
\eqref{E:mass-sum}, we necessarily have $M[\psi^j] < 1$ for each
$j$, which by \eqref{E:energy-below} implies that for $n$
sufficiently large,
$$
M[e^{-it_n^j\Delta}\psi^j] E[e^{-it_n^j\Delta}\psi^j] <
(ME)_\text{c}.
$$
For a given $j$, there are two cases to consider: Case (a). If
$|t_n^j| \to +\infty$ (passing to a subsequence we have $t_n^j\to
+\infty$ or $t_n^j\to -\infty$) we have
$\|e^{-it_n^j\Delta}\psi^j\|_{L^4} \to 0$ (as discussed in the proof
of Corollary \ref{C:energy_expand}), and thus,
$$
\frac12\|\psi^j\|_{L^2}^2\|\nabla\psi^j\|_{L^2}^2< (ME)_\text{c}
$$
(we have used $\|\nabla e^{-it_n^j\Delta}\psi^j\|_{L^2} = \|\nabla
\psi^j\|_{L^2}$). Let $\nls(t)\psi$ denote the solution to
\eqref{E:NLS} with initial data $\psi$.  By the existence of wave
operators (Proposition \ref{P:waveop}), there exists $\tilde \psi^j$
such that
$$
\|\nls(-t_n^j)\tilde \psi^j - e^{-it_n^j\Delta}\psi^j\|_{H^1} \to 0,
\quad \text{as }n\to +\infty
$$
with
$$
\|\tilde \psi^j\|_{L^2} \,\|\nabla \, \nls(t)\tilde\psi^j\|_{L^2} <
\|Q\|_{L^2}\|\nabla Q\|_{L^2},
$$
$$
M[\tilde \psi^j] = \|\psi^j\|_{L^2}^2, \quad E[\tilde\psi^j] =
\frac12 \, \|\nabla \psi^j\|_{L^2}^2,
$$
and thus,
$$
M[\tilde \psi^j] E[\tilde \psi^j] < (ME)_\text{c}, \quad
\|\nls(t)\tilde \psi^j\|_{S(\dot H^{1/2})} < +\infty.
$$
Case (b). On the other hand, if for a given $j$ we have $t_n^j\to
t_*$ finite (and there can be at most one such $j$ by
\eqref{E:pwdiverge}), we note that by continuity of the linear flow in $H^1$,
$$
e^{-it_n^j\Delta}\psi^j \to e^{-it_*\Delta}\psi^j \quad
\text{strongly in }H^1,
$$
and we let $\tilde \psi^j=\nls(t_*)[e^{-it_*\Delta}\psi^j]$ so that
$\nls(-t_*)\tilde \psi^j = e^{-it_*\Delta}\psi^j$.  In either
case, associated to each original profile $\psi^j$ we now have a new
profile $\tilde \psi^j$ such that
$$
\|\nls(-t_n^j)\tilde \psi^j - e^{-it_n^j\Delta}\psi^j \|_{H^1} \to 0
\quad \text{as }n\to +\infty.
$$
It now follows that we can replace $e^{-it_n^j\Delta}\psi^j$ by
$\nls(-t_n^j)\tilde \psi^j$ in \eqref{E:profexpand} to obtain
$$
u_{n,0} = \sum_{j=1}^M \nls(-t_n^j)\tilde \psi^j + \tilde W_n^M,
$$
where
$$
\lim_{M\to +\infty} \left[ \lim_{n\to +\infty} \|e^{it\Delta}\tilde
W_n^M\|_{S(\dot H^{1/2})} \right] =0.
$$
The idea of what follows is that we approximate
$$
\nls(t)u_{n,0} \approx \sum_{j=1}^M \nls(t-t_n^j)\tilde \psi^j
$$
via a perturbation theory argument, and since the right side has
bounded $S(\dot H^{1/2})$ norm, so must the left-side, which is a
contradiction.  To carry out this argument, we introduce the
notation $v^j(t)=\nls(t)\tilde \psi^j$, $u_n(t)=\nls(t)u_{n,0}$,
and\footnote{$\tilde u_n$, and $e_n$ also depend on $M$, but we have
suppressed the notation.}
$$
\tilde u_n(t) = \sum_{j=1}^M v^j(t-t_n^j).
$$
Then
$$
i\partial_t \tilde u_n + \Delta \tilde u_n + |\tilde u_n|^2 \tilde
u_n = e_n,
$$
where
$$
e_n = |\tilde u_n|^2 \tilde u_n - \sum_{j=1}^M |v^j(t-t_n^j)|^2
v^j(t-t_n^j).
$$
We claim that there is a (large) constant $A$ (independent of $M$)
with the property that for any $M$, there exists $n_0=n_0(M)$ such
that for $n>n_0$,
$$
\|\tilde u_n\|_{S(\dot H^{1/2})} \leq A .
$$
Moreover, we claim that for each $M$ and $\epsilon>0$ there exists
$n_1=n_1(M,\epsilon)$ such that for $n>n_1$,
$$
\|e_n\|_{L_t^{10/3}L_x^{5/4}} \leq \epsilon.
$$
Note that since $\tilde u_n(0)-u_n(0)=\tilde W_n^M$, there exists
$M_1=M_1(\epsilon)$ sufficiently large such that for each $M>M_1$
there exists $n_2=n_2(M)$ such that $n>n_2$ implies
$$
\|e^{it\Delta}(\tilde u_n(0)-u_n(0))\|_{S(\dot H^{1/2})} \leq
\epsilon.
$$
Thus, we may apply Proposition \ref{P:longtime} (long-time
perturbation theory) to obtain that for $n$ and $M$ sufficiently
large, $\|u_n\|_{S(\dot H^{1/2})}<\infty$, a
contradiction.\footnote{The order of logic here is:  The constant
$A$, which is independent of $M$, is put into Prop.
\ref{P:longtime}, which gives a suitable $\epsilon$.  We then take
$M_1=M_1(\epsilon)$ as above, and then take $n=\max(n_0,n_1,n_2)$.}

Therefore, it remains to establish the above claims, and we begin
with showing that $\|\tilde u_n\|_{S(\dot H^{1/2})} \leq A$  for
$n>n_0=n_0(M)$, where $A$ is some large constant independent of $M$.
Let $M_0$ be large enough so that
$$
\|e^{it\Delta}\tilde W_n^{M_0}\|_{S(\dot H^{1/2})} \leq
\delta_{\text{sd}}.
$$
Then for each $j>M_0$, we have $\|e^{it\Delta}\psi^j\|_{S(\dot
H^{1/2})} \leq \delta_{\text{sd}}$, and by the second part of
Proposition \ref{P:waveop} we obtain
\begin{equation}
 \label{E:twice}
\|v^j\|_{S(\dot H^{1/2})} \leq 2\|e^{it\Delta}\psi^j\|_{S(\dot
H^{1/2})} \quad \text{for }j> M_0.
\end{equation}
By the elementary inequality: for $a_j\geq 0$,
$$
\Big| \Big( \sum_{j=1}^M a_j \Big)^{5/2} - \sum_{j=1}^M a_j^{5/2}
\Big| \leq c_M \sum_{j\neq k} |a_j| |a_k|^{3/2} \, ,
$$
we have
\begin{equation}
\label{E:u_nbound}
\begin{aligned}
\|\tilde u_n\|_{L_t^5 L_x^5}^5
&= \sum_{j=1}^{M_0} \|v^j\|_{L_t^5 L_x^5}^5 + \sum_{j=M_0+1}^M \|v^j\|_{L_t^5 L_x^5}^5
+ \text{cross terms} \\
&\leq \sum_{j=1}^{M_0} \|v^j\|_{L_t^5 L_x^5}^5 + 2^5\sum_{j=M_0+1}^M
\|e^{it\Delta}\psi^j \|_{L_t^5 L_x^5}^5 + \text{cross terms} \,,
\end{aligned}
\end{equation}
where we used (\ref{E:twice}) to bound middle terms. On the other
hand, by \eqref{E:profexpand},
\begin{equation}
\label{E:linexpand} \|e^{it\Delta}u_{n,0}\|_{L_t^5 L_x^5}^5 =
\sum_{j=1}^{M_0} \|e^{it\Delta}\psi^j\|_{L_t^5 L_x^5}^5 +
\sum_{j=M_0+1}^M \|e^{it\Delta}\psi^j \|_{L_t^5 L_x^5}^5 +
\text{cross terms}\, .
\end{equation}
The ``cross terms'' are made $\leq 1$ by taking $n_0=n_0(M)$ large
enough and appealing to \eqref{E:pwdiverge}. We observe that since
$\|e^{it\Delta}u_{n,0}\|_{L_t^5L_x^5} \leq c\|u_{n,0}\|_{\dot
H^{1/2}} \leq c'$, \eqref{E:linexpand} shows that the quantity
$\sum_{j=M_0+1}^M \|e^{it\Delta}\psi^j \|_{L_t^5 L_x^5}^5$ is
bounded independently of $M$ provided $n>n_0$.  Then,
\eqref{E:u_nbound} gives that $\|\tilde u_n\|_{L_t^5 L_x^5}$ is
bounded independently of $M$ provided $n>n_0$. A similar argument
establishes that $\|\tilde u_n\|_{L_t^\infty L_x^3}$ is bounded
independently of $M$ for $n>n_0$.  Interpolation between these
exponents gives that $\|\tilde u_n\|_{L_t^{20}L_x^{10/3}}$  is
bounded independently of $M$ for $n>n_0$.  Finally, by applying the
Kato estimate \eqref{E:Kato} to the integral equation for
$i\partial_t \tilde u_n + \Delta \tilde u_n + |\tilde u_n|^2 \tilde
u_n =e_n$ and using that $\|e_n\|_{S(\dot H^{-1/2})}\leq 1$ (proved
next), we obtain that $\|\tilde u_n\|_{S(\dot H^{1/2})}$  is bounded
independently of $M$ for $n>n_0$.

We now address the next claim, that for each $M$ and $\epsilon>0$, there
exists $n_1=n_1(M,\epsilon)$ such that for
$n>n_1$, $\|e_n\|_{L_t^{10/3}L_x^{5/4}} \leq \epsilon$.  The
expansion of $e_n$ consists of $\sim M^3$ cross terms of the form
$$
v^j(t-t^j_n) v^k(t-t_n^k) v^\ell(t-t_n^\ell),
$$
where not all three of $j$, $k$, and $\ell$ are the same.  Assume,
without loss, that $j\neq k$, and thus, $|t_n^j-t_n^k| \to \infty$
as $n\to +\infty$. We estimate
$$
\| v^j(t-t^j_n) v^k(t-t_n^k)
v^\ell(t-t_n^\ell)\|_{L_t^{10/3}L_x^{5/4}} \leq \| v^j(t-t^j_n)
v^k(t-t_n^k) \|_{L_t^{10}L_x^{5/3}}
\|v^\ell(t-t_n^\ell)\|_{L_t^5L_x^5}.
$$
Now observe that
$$
\|v^j(t-(t_n^j-t_n^k))\cdot v^k(t)\|_{L_t^{10}L_x^{5/3}}\to 0,
$$
since $v^j$ and $v^k$ belong to $L_t^{20}L_x^{10/3}$ and
$|t_n^j-t_n^k| \to \infty$.

\textit{Case 2}.  $\psi^1\neq 0$, and $\psi^j=0$ for all $j\geq 2$.

By \eqref{E:mass-sum}, we have $M[\psi^1] \leq 1$ and by
\eqref{E:energy-below}, we have $\lim_{n\to +\infty}
E[e^{-it_n^1\Delta}\psi^1] \leq (ME)_\crit$.  If $t_n^1$ converges
(to $0$ without loss of generality), we take $\tilde \psi^1 =
\psi^1$ and then we have $\|\nls(-t_n^1)\tilde \psi^1 -
e^{-it_n^1\Delta}\psi^1\|_{H^1} \to 0$ as $n\to +\infty$.  If, on
the other hand, $t_n^1\to +\infty$, then since
$\|e^{it_n^1\Delta}\psi^1\|_{L^4} \to 0$,
$$
\frac12\|\nabla \psi^1 \|_{L^2}^2 = \lim_{n\to+\infty}
E[e^{-it_n^1\Delta}\psi^1] \leq (ME)_\crit.
$$
Thus, by the existence of wave operators (Proposition
\ref{P:waveop}), there exists $\tilde \psi^1$ such that $M[\tilde
\psi^1] = M[\psi^1] \leq 1$,  $E[\tilde \psi^1] =
\frac12\|\nabla \psi^1\|_{L^2}^2 \leq (ME)_\crit$, and
$\|\nls(-t_n^1)\tilde \psi^1 - e^{-it_n^1\Delta}\psi^1\|_{H^1} \to
0$ as $n\to +\infty$.

In either case, let $\tilde W_n^M = W_n^M + (e^{-it_n^1\Delta}\psi^1
- \nls(-t_n^1)\tilde \psi^1)$. Then, by the Strichartz estimates,
$$
\|e^{-it\Delta}\tilde W_n^M\|_{S(\dot H^{1/2})} \leq
\|e^{-it\Delta}W_n^M\|_{S(\dot H^{1/2})} +
c\|e^{-it_n^1\Delta}\psi^1 - \nls(-t_n^1)\tilde \psi^1\|_{\dot
H^{1/2}}\, ,
$$
and therefore, $\lim_{n\to +\infty} \|e^{-it\Delta}\tilde
W_n^M\|_{S(\dot H^{1/2})} = \lim_{n\to +\infty}
\|e^{-it\Delta}W_n^M\|_{S(\dot H^{1/2})}$. Hence, we now have
$$
u_{n,0} = \nls(-t_n^1)\tilde\psi^1 + \tilde W_n^M
$$
with $M[\tilde \psi^1] \leq 1$, $E[\tilde \psi^1] \leq (ME)_\crit$, and
$$
\limsup_{M\to +\infty}\Big[ \lim_{n\to +\infty}
\|\tilde W_n^M\|_{S(\dot H^{1/2})}\Big] = 0 \, .
$$

Let $u_\crit$ be the solution to \eqref{E:NLS} with initial data
$u_{\crit,0}=\tilde\psi^1$. Now we claim that $\|u_\crit\|_{S(\dot
H^{1/2})} = \infty$, and thus, $M[u_\crit]=1$ and
$E[u_\crit]=(ME)_\crit$, which will complete the proof.  To
establish this claim, we use a perturbation argument similar to that
in Case 1.  Suppose
$$
A := \|\nls(t-t_n^1)\tilde\psi^1\|_{S(\dot H^{1/2})}
= \|\nls(t)\tilde\psi^1\|_{S(\dot H^{1/2})}
=\|u_\crit\|_{S(\dot H^{1/2})} < \infty.
$$
Obtain $\epsilon_0=\epsilon_0(A)$ from the long-time perturbation
theory (Proposition \ref{P:longtime}), and then take $M$
sufficiently large and $n_2=n_2(M)$ sufficiently large so that
$n>n_2$ implies $\|\tilde W_n^M\|_{S(\dot H^{1/2})} \leq
\epsilon_0$.  We then repeat the argument in Case 1 using
Proposition \ref{P:longtime} to obtain that there exists $n$ large
for which $\|u_n\|_{S(\dot H^{1/2})} <\infty$, a contradiction.

\end{proof}

\begin{proposition}[Precompactness of the flow of the critical solution]
\label{P:crit_compact}
With $u_\crit$ as in Proposition \ref{P:exist_crit}, let
$$
K = \{ \, u_\crit(t) \, | \, t\in [0,+\infty) \, \} \subset H^1.
$$
Then $K$ is precompact in $H^1$ (i.e.\ $\bar K$ is compact in $H^1$).
\end{proposition}
\begin{proof}
Take a sequence $t_n \to
+\infty$; we shall argue that $u_\crit(t_n)$ has a subsequence
converging in $H^1$.\footnote{By time continuity of the solution in
$H^1$, we of course do not need to consider the case when $t_n$ is
bounded and thus has a subsequence convergent to some finite time.}
Take $\phi_n = u_\crit(t_n)$ (a uniformly bounded sequence in $H^1$)
in the profile expansion lemma (Lemma \ref{L:pe}) to obtain profiles
$\psi^j$ and an error $W_n^M$ such that
$$
u_\crit(t_n) = \sum_{j=1}^M e^{-it_n^j\Delta}\psi^j + W_n^M
$$
with $|t_n^j-t_n^k| \to +\infty$ as $n\to +\infty$ for fixed $j\neq
k$. By the energy Pythagorean expansion (Corollary
\ref{C:energy_expand}), we have
$$
\sum_{j=1}^M \lim_{n\to +\infty} E[e^{-it_n^j\Delta}\psi^j] +
\lim_{n\to +\infty} E[W_n^M] = E[u_\crit] = (ME)_\crit,
$$
and thus (recalling that each energy is $\geq 0$ -- see Lemma
\ref{L:comparability}),
$$
\lim_{n\to +\infty} E[e^{-it_n^j\Delta}\psi^j] \leq (ME)_\crit \quad
\forall \; j \,.
$$
Also by $s=0$ of \eqref{E:pythag}, we have
$$
\sum_{j=1}^M M[\psi^j] + \lim_{n\to +\infty} M[W_n^M] = \lim_{n\to
+\infty} M[u_{n,0}] =1.
$$
We now consider two cases, just as in the proof of Proposition
\ref{P:exist_crit}; both Case 1 and Case 2 will lead to a
contradiction.

\textit{Case 1}. More than one $\psi^j \neq 0$.  The proof that this
leads to a contradiction is identical to that in Proposition
\ref{P:exist_crit}, so we omit it.

\textit{Case 2}. Only $\psi^1 \neq 0$ and $\psi^j =0$ for all $2
\leq j \leq M$, so that
\begin{equation}
\label{E:case2expand}
u_\crit(t_n) = e^{-it_n^1\Delta}\psi^1 +W_n^M
\end{equation}
Just as in the proof of Proposition \ref{P:exist_crit} Case 2, we obtain that
$$
M[\psi^1] =1, \quad \lim_{n\to +\infty} E[e^{-it_n^1\Delta}\psi^1] =
(ME)_\crit \,,
$$
$$
\lim_{n\to +\infty} M[W_n^M] =0, \quad \text{and} \quad \lim_{n\to
+\infty} E[W_n^M] = 0 \,.
$$
By the comparability lemma (Lemma \ref{L:comparability}),
\begin{equation}
\label{E:H^1tozero}
\lim_{n\to +\infty} \| W_n^M\|_{H^1} =0 \,.
\end{equation}
Next, we show that (a subsequence of) $t_n^1$ converges.\footnote{In
the rest of the argument, take care not to confuse $t_n^1$
(associated with $\psi^1$) with $t_n$.} Suppose that $t_n^1 \to
-\infty$. Then
$$
\|e^{it\Delta}u_\crit(t_n)\|_{S(\dot H^{1/2}; [0,+\infty))} \leq \|
e^{i(t-t_n^1)\Delta}\psi^1\|_{S(\dot H^{1/2}; [0,+\infty))} + \|
e^{i t \Delta} W_n^M\|_{S(\dot H^{1/2}; [0,+\infty))} .
$$
Since
$$
\lim_{n\to +\infty} \| e^{i(t-t_n^1)\Delta}\psi^1\|_{S(\dot H^{1/2};
[0,+\infty))} = \lim_{n\to +\infty} \| e^{it\Delta}\psi^1\|_{S(\dot
H^{1/2}; [-t_n^1,+\infty))} = 0
$$
and $\|e^{it\Delta}W_n^M\|_{S(\dot H^{1/2})} \leq \frac12 \,
\delta_{\text{sd}}$, we obtain a contradiction to the small data
scattering theory (Proposition \ref{P:sd}) by taking $n$
sufficiently large.  On the other hand, suppose that $t_n^1 \to
+\infty$.  Then we can similarly argue that for $n$ large,
$$
\|e^{it\Delta}u_\crit(t_n)\|_{S(\dot H^{1/2}; (-\infty,0])} \leq
\frac12\delta_{\text{sd}},
$$
and thus, the small data scattering theory (Proposition \ref{P:sd})
shows that
$$
\|u_\crit\|_{S(\dot H^{1/2};(-\infty,t_n])} \leq
\delta_{\text{sd}}.
$$
Since $t_n\to +\infty$, by sending $n\to +\infty$ in the above, we
obtain $\|u_\crit\|_{S(\dot H^{1/2};(-\infty,+\infty))} \leq
\delta_{\text{sd}}$, a contradiction. Thus, we have shown that
$t_n^1$ converges to some finite $t^1$.

Since $e^{-it_n^1\Delta}\psi^1 \to e^{-it^1\Delta}\psi^1$ in $H^1$
and \eqref{E:H^1tozero} holds, \eqref{E:case2expand} shows that
$u_\crit(t_n)$ converges in $H^1$.
\end{proof}

\begin{lemma}[Precompactness of the flow implies uniform localization]
\label{L:unif_small} Let $u$ be a solution to \eqref{E:NLS} such that
$$K = \{ \, u(t) \, | \, t\in [0,+\infty) \, \}$$
is precompact in $H^1$.  Then for each $\epsilon>0$, there exists $R> 0$
so that
$$
\int_{|x|>R} |\nabla u(x,t)|^2 \leq \epsilon, \quad \text{for
all }0\leq t <+\infty.
$$
\end{lemma}
\begin{proof}
If not, then there exists $\epsilon>0$ and a sequence of times $t_n$
such that
$$
\int_{|x|>n} |\nabla u(x,t_n)|^2 \, dx \geq \epsilon .
$$
Since $K$ is precompact, there exists $\phi\in H^1$ such that,
passing to a subsequence of $t_n$, we have  $u(t_n) \to \phi$ in
$H^1$.  By taking $n$ large, we have both
$$
\int_{|x|>n} |\nabla \phi(x)|^2 \leq \frac14\,\epsilon
$$
and
$$
\int_{\mathbb{R}^3} |\nabla(u(x,t_n)-\phi(x))|^2 \, dx \leq
\frac14\,\epsilon\,,
$$
which is a contradiction.
\end{proof}

\section{Rigidity theorem}
\label{S:rigidity}

We now prove the rigidity theorem.
\begin{theorem}[Rigidity]
\label{T:rigidity}
Suppose $u_0 \in H^1$ satisfies
\begin{equation}
 \label{E:comp1}
M[u_0]E[u_0] < M[Q]E[Q]
\end{equation}
and
\begin{equation}
 \label{E:comp2}
\|u_0\|_{L^2} \|\nabla u_0\|_{L^2} < \|Q\|_{L^2} \|\nabla Q
\|_{L^2}.
\end{equation}
Let $u$ be the global $H^1$ solution of \eqref{E:NLS} with initial data $u_0$ and
suppose that
$$
K = \{ \, u(t) \, | \, t \in [0,+\infty) \, \} \quad \text{is
precompact in} ~ H^1.
$$
Then $u_0 = 0$.
\end{theorem}

\begin{proof}
Let $\phi \in C_0^{\infty}$, radial, with
$$ \phi(x)=\left\{
\begin{array}{lll}
|x|^2 &\text{for} &|x|\leq 1\\
0 &\text{for}&|x|\geq 2
\end{array}
\right. \,.$$
For $R>0$ define $z_R(t) = \int R^2 \phi(\frac{x}{R}) \,
|u(x,t)|^2 \, dx$. Then
\begin{equation}
 \label{E:localvirial1}
|z'_R(t)| \leq 2 R \left| \int \bar{u}(t)\, \nabla u(t) \, (\nabla
\phi)\left(\frac{x}{R}\right) \, dx \right| \leq c \, R \int_{0 <
|x|< 2R} |\nabla u(t)| \, |u(t)| \, dx.
\end{equation}
Using H\"older's inequality and Theorem
\ref{T:main}(\ref{I:global})(a), we bound the previous expression
by
$$
c \, R \, \|\nabla u(t)\|_{L^2} \, \|u \|_{L^2} \leq c \, R \,
\|\nabla Q\|_{L^2} \, \|Q\|_{L^2} = \tilde c \, R.
$$
Thus, we obtain
\begin{equation}
 \label{E:firstderivative}
|z'_R(t) - z'_R(0)| \leq 2 \, \tilde c \, R \quad \text{for} \quad t
> 0.
\end{equation}

Next we estimate $z''_R(t)$ using the localized virial identity
\eqref{E:localvirial}:
\begin{align*}
z''_R(t)
& = 4 \int \phi''\left(\frac{|x|}{R}\right) \,|\nabla u|^2
- \frac1{R^2} \int (\Delta^2 \phi) \left(\frac{x}{R}\right) |u|^2 -
\int (\Delta \phi)\left(\frac{x}{R}\right) \,
|u|^{4}\\
& \geq 8 \int_{|x|\leq R} |\grad u|^2 + 4 \int_{R<|x|<2R}
\phi''\left(\frac{|x|}{R}\right) \,|\nabla u|^2 -   \frac{c}{R^2}
\int_{R<|x|<2R} {|u|^2}\\
&\qquad - 6 %
\int_{|x| \leq R} |u|^{4} - c\, \int_{R<|x|<2R} |u|^{4}\\
& \geq \left( 8 \int_{|x| \leq R} |\nabla u|^2 - 6 \int_{|x| \leq R}
|u|^{4} \right) - c_1 \, \int_{R<|x|<2R} \left(|\nabla u|^2 +
\frac{|u|^2}{R^2} + |u|^{4} \right).
\end{align*}

Since \eqref{E:comp1} holds, take $\delta>0$ such that $M[u_0]E[u_0]
\leq (1-\delta) M[Q]E[Q]$. Let $\epsilon = c_1^{-1} \, c_\delta \,
\int |\nabla u_0|^2$, where $c_\delta$ is as in \eqref{E:comp3}.

Since $\{u(t) | t \in [0, \infty)\}$ is precompact in $H^1$, by
Lemma \ref{L:unif_small} there exists $R_1>0$ such that
$\int_{|x|>R_1} |\nabla u(t)|^2 \leq \frac19 \,\epsilon$. Next,
because of mass conservation, there exists $R_2> 0$ such that
$\frac1{R^2_2} \int |u|^2 < \frac19 \, \epsilon$. Finally, the
radial Gagliardo-Nirenberg inequality \eqref{E:Strauss} yields the
existence of $R_3>0$ such that
$$
\int_{|x|>R_3} |u(t)|^4 \leq \frac{c}{R^2_3} \, \|\nabla
u(t)\|_{L^2(|x|>R_3)} \, \|u_0\|_{L^2}^3 \leq \frac{c}{R^2_3} \,
\|\nabla u_0 \|_{L^2} \, \|u_0\|_{L^2}^3 \leq \frac19 \,\epsilon,
$$
with $R_3^2 > 9 \, c\, \ep \, \|\nabla u_0 \|_{L^2} \,
\|u_0\|_{L^2}^3$; in the above chain we used the gradient-energy
comparability (Lemma \ref{L:comparability}) with $t=0$ on the left
side. Take $R = \max \{R_1, R_2, R_3\}$ to obtain
\begin{equation}
 \label{E:allterms}
c_1 \, \int_{|x| > R} \left(|\nabla u|^2 + \frac{|u|^2}{R^2} +
|u|^{4} \right) \leq \frac13 \, c_\delta \int |\nabla u_0|^2.
\end{equation}
By \eqref{E:comp3} and Lemma
\ref{L:comparability}, we also have
\begin{equation}
 \label{E:bigterms}
8 \int|\nabla u|^2 - 6 \int |u|^{4} \geq
 c_\delta \, \int |\nabla u_0|^2.
\end{equation}
Splitting the integrals on the left side of the above expression
into the regions $\{|x|>R\}$ and $\{|x|<R\}$ and applying
(\ref{E:allterms}), we get
$$
8 \int_{|x| \leq R} |\nabla u|^2 - 6 \int_{|x| \leq R} |u|^{4} \geq
 \frac23 \, c_\delta \, \int |\nabla u_0|^2 .
$$
Hence, we obtain $ z''_R(t) \geq \frac13 \, c_\delta \, \|\nabla u_0
\|_{L^2}^{2}$, which implies by integration from 0 to $t$ that
$z'_R(t) - z'_R(0) \geq \frac13 \, c_\delta \, \|\nabla u_0
\|_{L^2}^{2} \, t$. Taking $t$ large, we obtain a contradiction with
(\ref{E:firstderivative}), which can be resolved only if $\|\nabla
u_0 \|_{L^2} = 0.$
\end{proof}

To complete the proof of Theorem \ref{T:main}\eqref{I:global}(b), we
just apply Theorem \ref{T:rigidity} to $u_\crit$ constructed in
Proposition \ref{P:exist_crit}, which by Proposition
\ref{P:crit_compact}, meets the hypotheses in Theorem
\ref{T:rigidity}.  Thus $u_{\crit,0}=0$, which contradicts the fact
that $\|u_\crit\|_{S(\dot H^{1/2})}=\infty$. We have thus obtained
that if $\|u_0\|_{L^2}\|\nabla u_0\|_{L^2}<\|Q\|_{L^2}\|\nabla
Q\|_{L^2}$ and $M[u]E[u]<M[Q]E[Q]$, then $\text{SC}(u_0)$ holds,
i.e.\ $\|u\|_{S(\dot H^{1/2})} <\infty$.  By Proposition
\ref{P:persistence}, $H^1$ scattering holds.

\section{Extensions to general mass supercritical, energy subcritical NLS equations}
\label{S:extensions}

Consider the focusing mass supercritical, energy subcritical
nonlinear Schr\"odinger equation NLS$_p(\cR^N)$:
\begin{equation}
\label{E:NLSgeneral}
\left\{ \begin{array}{l} i\partial_t u +\Delta u + |u|^{p-1}u=0,
\quad (x,t) \in \cR^N
\times \cR,\\
u(x,0) = u_0(x) \in H^1(\cR^N),
\end{array}
\right.
\end{equation}
with the choice of nonlinear exponent $p$ and the dimension $N$ such that
$$
0<s_c<1, \quad \text{where} \quad s_c = \frac{N}2-\frac2{p-1}.
$$
The initial value problem with $u_0 \in {H}^1(\cR^N)$ is locally
well-posed, see \cite{GV79}. Denote by $I = (-T_*, T^*)$ the maximal
interval of existence of the solution $u$ (e.g., see
\cite{Caz-book}). This implies that either $T^* = +\infty$ or $T^* <
+\infty$ and $\Vert \nabla u (t) \Vert_{L^2} \to \infty$ as $t \to
T^*$ (similar properties for $T_*$).

The solutions to this problem satisfy mass and energy conservation
laws, in particular,
$$
E[u(t)] = \frac12\int |\nabla u(x,t)|^2 - \frac1{p+1}\int
|u(x,t)|^{p+1} \,dx = E[u_0].
$$
The Sobolev $\dot{H}^{s_c}$ norm
is invariant under the scaling $u \mapsto u_{\lambda}(x,t) =
\lambda^{2/(p-1)} u(\lambda x, \lambda^2 t)$ ($u_\lambda$ is a
solution of NLS$_p(\cR^N)$, if $u$ is).

The general Gagliardo-Nirenberg inequality (see \cite{W83}) is valid
for values of $p$ and $N$ such that $0 \leq s_c < 1$\footnote{It is
also valid for $s_c=1$ becoming the Sobolev embedding, see Remark
\ref{W<->KM}.}:
\begin{equation}
 \label{GNgeneral}
\Vert u \Vert^{p+1}_{L^{p+1}(\cR^N)} \leq c_{\text{GN}} \, \Vert \nabla u
\Vert_{L^2(\cR^N)}^{\frac{N(p-1)}{2}} \, \Vert u
\Vert_{L^2(\cR^N)}^{2-\frac{(N-2)(p-1)}{2}},
\end{equation}
where
$$
c_{\text{GN}} = \frac{\Vert Q \Vert^{p+1}_{L^{p+1}(\cR^N)}}{\Vert
\nabla Q \Vert_{L^2(\cR^N)}^{\frac{N(p-1)}{2}} \, \Vert Q
\Vert_{L^2(\cR^N)}^{2-\frac{(N-2)(p-1)}{2}}}
$$
and $Q$ is the ground state solution (positive solution of minimal
$L^2$ norm) of the equation
\begin{equation}
 \label{E:Qground}
-(1-s_c)Q+\Delta Q + \,|Q|^{p-1} \, Q = 0.
\end{equation}
(See \cite{W83} and references therein for discussion on the
existence of positive solutions of class $H^1(\cR^N)$ to this
equation.)\footnote{In the case $p=3$, $N=3$, we have $s_c=\frac12$,
and thus, the normalization for $Q$ chosen here is different from
that in the main part of this paper.  The normalization of $Q$ taken
here was chosen since it enables us to draw a comparison with the
$s_c=1$ endpoint result of Kenig-Merle \cite{KM06a}.}  The
corresponding soliton solution to \eqref{E:NLSgeneral} is $u(x,t) =
e^{i(1-s_c)t}Q(x)$.

The generalization of Theorem \ref{BlowUp} (or Theorem \ref{T:main}
without scattering) to all $0< s_c < 1$ is based on using the
scaling invariant quantity $\Vert \nabla u \Vert^{s_c}_{L^2(\cR^N)}
\cdot \Vert u \Vert_{L^2(\cR^N)}^{1-s_c}$.

\begin{theorem}
 \label{T:BlowUpN}
Consider NLS$_p(\cR^N)$ with (possibly non-radial) $u_0 \in
{H}^1(\cR^N)$  and $0 < s_c < 1$. Suppose that
\begin{equation}
 \label{E:initialLambda}
\qquad  E[u_0]^{s_c}\, M[u_0]^{1-s_c} < E[Q]^{s_c} \, M[Q]^{1-s_c}, \quad E[u_0] \geq 0.\\
\end{equation}
If \eqref{E:initialLambda} holds and
\begin{equation}
 \label{E:less}
\Vert \nabla u_0 \Vert^{s_c}_{L^2(\cR^N)} \Vert u_0
\Vert_{L^2(\cR^N)}^{1-s_c} < \Vert \nabla Q \Vert^{s_c}_{L^2(\cR^N)}
\Vert Q \Vert_{L^2(\cR^N)}^{1-s_c},
\end{equation}
then for any $t\in I$,
\begin{equation}
 \label{E:lessallt}
\Vert \nabla u(t) \Vert^{s_c}_{L^2(\cR^N)} \Vert u_0
\Vert_{L^2(\cR^N)}^{1-s_c} < \Vert \nabla Q \Vert^{s_c}_{L^2(\cR^N)}
\Vert Q \Vert_{L^2(\cR^N)}^{1-s_c},
\end{equation}
and thus $I = (-\infty, +\infty)$, i.e. the solution exists globally
in time.\\
If \eqref{E:initialLambda} holds and
\begin{equation}
 \label{E:greater}
\Vert \nabla u_0 \Vert^{s_c}_{L^2(\cR^N)} \Vert u_0
\Vert_{L^2(\cR^N)}^{1-s_c} > \Vert \nabla Q \Vert^{s_c}_{L^2(\cR^N)}
\Vert Q \Vert_{L^2(\cR^N)}^{1-s_c},
\end{equation}
then for $t \in I$
\begin{equation}
 \label{E:greaterallt}
\Vert \nabla u(t) \Vert^{s_c}_{L^2(\cR^N)} \Vert u_0
\Vert^{1-s_c}_{L^2(\cR^N)} > \Vert \nabla Q \Vert^{s_c}_{L^2(\cR^N)}
\Vert Q \Vert_{L^2(\cR^N)}^{1-s_c}.
\end{equation}
Furthermore, if (a) $|x|u_0 \in L^2(\cR^N)$, or (b) $u_0$ is radial
with $N > 1$ and $1+\frac4{N} < p < \min\{1+ \frac4{N-2}, 5 \}$,
then $I$ is finite, and thus, the solution blows up in finite time.
The finite-time blowup conclusion and \eqref{E:greaterallt} also
hold if, in place of \eqref{E:initialLambda} and \eqref{E:greater},
we assume $E[u_0]<0$.
\end{theorem}

The proof of this theorem is similar to Theorem \ref{BlowUp} and can
be found in \cite{HR06}.

\begin{remark}
A finite-time $T$ blow-up solution to a mass-supercritical energy
subcritical NLS equation satisfies a lower bound on the blow-up
rate: $\|\nabla u(t)\|_{L^2} \geq c(T-t)^{-\alpha}$, where
$\alpha=\alpha(p,d)$.  This is obtained by scaling the local-theory,
and it implies that the quantity $\|u_0\|_{L^2}^{s_c}\|\nabla
u(t)\|_{L^2}^{1-s_c} \to \infty$, thus strengthening the conclusion
\eqref{E:greaterallt}.  A stronger result in this direction was
recently obtained by Merle-Rapha\"el \cite{MR06}:  if $u(t)$
blows-up in finite time $T^*<\infty$, then $\lim_{t\to T^*}
\|u(t)\|_{\dot H^{s_c}} = \infty$ (in fact, it diverges to $\infty$
with a logarithmic lower bound).
\end{remark}

\begin{remark} \label{W<->KM}
This theorem provides a link between the mass critical NLS and
energy critical NLS equations: Consider $s_c=1$; the theorem holds
true by the work of Kenig-Merle \cite[Section 3]{KM06a}. In this
case there is no mass involved, the Gagliardo-Nirenberg inequality
(\ref{GNgeneral}) becomes the Sobolev inequality, the condition
(\ref{E:initialLambda}) is $E[u_0] < E[Q]$, where $Q$ is the radial
positive decreasing (class $\dot{H}^1(\cR^N)$) solution of $\lap Q +
|Q|^{p-1} Q = 0$, and the conditions (\ref{E:less}) -
(\ref{E:greaterallt}) involve only the size of $\Vert \nabla u_0
\Vert_{L^2}$ in relation to $\Vert \nabla Q \Vert_{L^2(\cR^N)}$.  In
regard to the case $s_c=0$, \eqref{E:initialLambda} should be
replaced by $M[u]<M[Q]$ and \eqref{E:less} becomes the same
statement.  Under this hypotheses, the result of M.\ Weinstein
\cite{W83} states that
$$
\|\nabla u(t)\|_{L^2}^2 \leq 2\,\Big( 1 -
\frac{\|u_0\|_{L^2}^2}{\|Q\|_{L^2}^2} \Big)^{-1}E[u], \qquad
E[u]>0,
$$
and thus, global existence holds.  We do not recover this estimate
as a formal limit in \eqref{E:lessallt},\footnote{It might appear as
a formal limit if one were to refine the estimate \eqref{E:lessallt}
to account for the gain resulting from the strict inequality in
\eqref{E:initialLambda} (as we did in the proof of Theorem
\ref{BlowUp}) before passing to the $s_c\to 0$ limit.} however, the
conclusion about the global existence in this case does hold true.
Our intention here is not to \textit{reprove} the $s_c=0$ endpoint
result -- only to draw a connection to it. The hypothesis
\eqref{E:greater} should be replaced by its formal limit
$\|u_0\|_{L^2} > \|Q\|_{L^2}$, which is the complement of
\eqref{E:initialLambda}. Thus, the only surviving claim in Theorem
\ref{T:BlowUpN} regarding blow-up in the $s_c=0$ limit is that it
should hold under the hypothesis $E[u_0]<0$. Blow-up under this
hypothesis is the classical result of Glassey \cite{G77} in the case
of finite variance, and in the radial case it is the result of
Ogawa-Tsutsumi \cite{OT91}.
\end{remark}

We expect that the proof of scattering for NLS$_p(\cR^N)$ with $u_0
\in H^1(\cR^N)$ and $0 < s_c < 1$ when \eqref{E:initialLambda} and
\eqref{E:less} hold will carry over analogously to the $N=3$, $p=3$
case, provided (i) $N>1$ (the radial assumption in 1D does not help
to eliminate the translation defect of compactness); (ii) the Kato
estimate (as in \eqref{E:Kato}) or the more refined Strichartz
estimates by Foschi \cite{F04} are sufficient to complete the long
term perturbation argument\footnote{It may be necessary, for
example, to express the estimates in terms of the norm $\|
D^{s_c-\alpha} (\cdots) \|_{S(\dot H^\alpha)}$ for some $0<\alpha <
s_c$, rather than $\|\cdot \|_{S(\dot H^{s_c})}$.}.


\begin{thebibliography}{00}

\bibitem{BAK}
L. Berg\'e, T. Alexander, and Y. Kivshar,
\textit{Stability criterion for attractive Bose-Einstein condensates},
Phys. Rev. A, 62 (2000).

\bibitem{Begout}
P. B\'egout, \emph{Necessary conditions and sufficient conditions
for global existence in the nonlinear Schr\"odinger equation}, Adv.
Math. Sci. Appl. 12 (2002), no. 2, pp. 817--827.

\bibitem{Caz-book}
T.\ Cazenave,
\textit{Semilinear Schr\"odinger equations}.
Courant Lecture Notes in Mathematics, 10. New York University,
Courant Institute of Mathematical Sciences, New York;
American Mathematical Society, Providence, RI, 2003. xiv+323 pp. ISBN: 0-8218-3399-5.

\bibitem{CKSTT04}
J.\ Colliander, M.\ Keel, G.\ Staffilani, H.\ Takaoka, T.\ Tao,
\textit{Global existence and scattering for rough solutions of a nonlinear
Schr\"odinger equation on $\Bbb R\sp 3$},
Comm. Pure Appl. Math. 57 (2004), no. 8, pp.\ 987--1014.

\bibitem{CKSTTAnnals}
J.\ Colliander, M.\ Keel, G.\ Staffilani, H.\ Takaoka, T.\ Tao,
\textit{Global well-posedness and scattering for the energy-critical nonlinear
Schr\"odinger equation in $R^3$}, arxiv.org preprint \texttt{arXiv:math/0402129v7 [math.AP]}.


\bibitem{JILA}
E. Donley, N. Claussen, S. Cornish, J. Roberts, E. Cornell, and C. Wieman,
\textit{Dynamics of collapsing and exploding Bose-Einstein condensates},
Nature 412 (2001) pp. 295.

\bibitem{F}
G. Fibich, \textit{Some modern aspects of self-focusing theory}, in
\textit{Self-Focusing: Past and Present}, R.W. Boyd, S.G. Lukishova,
Y.R. Shen, editors, to be published by Springer. Available at
\texttt{http://www.math.tau.ac.il/$\sim$fibich/publications.html}.

\bibitem{F04}
D.\ Foschi, \textit{Inhomogeneous Strichartz estimates},  J. Hyper.
Diff. Eq.  2  (2005),  no. 1, 1--24.

\bibitem{GV79} J.\ Ginibre and G.\ Velo,
\textit{On a class of nonlinear Schr\"odinger equation. I. The
Cauchy problems; II. Scattering theory, general case}, J. Func.
Anal. 32 (1979), 1-32, pp.\ 33-71.

\bibitem{GV85}  J.\ Ginibre and G.\ Velo,
\textit{Scattering theory in the energy space for a class of nonlinear Schr\"odinger equations},
J. Math. Pures Appl. (9) 64 (1985), no. 4, pp.\ 363--401.

\bibitem{G77} R.\ T.\ Glassey, \textit{On the blowing up of solutions
to the Cauchy problem for nonlinear Schr\"odinger equation},
J. Math. Phys., 18, 1977, 9, pp.\ 1794--1797.

\bibitem{HR06}
J.\ Holmer and S.\ Roudenko,
\textit{On blow-up solutions to the 3D
cubic nonlinear Schr\"odinger equation}, AMRX Appl. Math. Res.
Express, vol. 2007, article ID abm004, {\tt
doi:10.1093/amrx/abm004}.

\bibitem{Kato}  T.\ Kato,
\textit{An $L\sp {q,r}$-theory for nonlinear Schr\"odinger equations},
Spectral and scattering theory and applications,  pp.\ 223--238,
Adv. Stud. Pure Math., 23, Math. Soc. Japan, Tokyo, 1994.

\bibitem{K01}
S.\ Keraani, \textit{On the defect of compactness for the Strichartz
estimates of the Schr\"odinger equation}, J. Diff. Eq. 175 (2001),
pp.\ 353--392

\bibitem{KT98}
M.\ Keel and T.\ Tao, \textit{Endpoint Strichartz estimates}, Amer.
J. Math., 120 (1998), pp.\ 955--980.

\bibitem{KM06a}
C.E.\ Kenig and F.\ Merle,
\textit{Global well-posedness, scattering, and blow-up for the energy-critical
focusing nonlinear Schr\"odinger equation in the radial case},
Invent. Math. 166 (2006), no. 3, pp.\ 645--675.

\bibitem{KPV93}
C.\ Kenig, G.\ Ponce and L.\ Vega, \textit{Well-posedness and
scattering results for the generalized Korteweg-de Vries equation
via the contraction principle}, Comm. Pure Appl. Math.  46  (1993),
no. 4, 527--620.


\bibitem{26}
N.E. Kosmatova, V.F. Shvets and V.E. Zakharov, \textit{Computer
simulation of wave collapses in the nonlinear Schr\"odinger
equation}, Physica D 52 (1991), pp.\ 16--35.

\bibitem{KRRT}
E.A. Kuznetsov, J. Juul Rasmussen, K. Rypdal, S. K.  Turitsyn,
\emph{Sharper criteria for the wave collapse},
Phys. D 87 (1995), no. 1-4, pp. 273--284.

\bibitem{LP-book} F.\ Linares, G.\ Ponce,
\textit{Introduction to nonlinear dispersive equations},
Rio de Janeiro: IMPA, 2004, 243 pp.

\bibitem{MR06} F. Merle, P. Rapha\"el,
\textit{Blow-up of the critical norm for some radial $L^2$
supercritical nonlinear Schr\"odinger equations}, arxiv.org preprint
\texttt{arXiv:math/0605378v2 [math.AP]}.


\bibitem{OT91} T. Ogawa and Y. Tsutsumi,
\textit{Blow-Up of $H^1$ solution for the Nonlinear Schr\"odinger
Equation}, J. Diff. Eq. 92 (1991), pp.\ 317-330.

\bibitem{Sch}
Schlein, B. \textit{Derivation of the Gross-Pitaevskii hierarchy},
Mathematical physics of quantum mechanics, pp. 279--293, Lecture
Notes in Phys., 690, Springer, Berlin, 2006. See also arxiv.org {\tt
arXiv:math-ph/0504078}.

\bibitem{Soffer}
A. Soffer, \emph{Soliton dynamics and scattering}, Proceedings of the
International Congress of Mathematicians, Madrid, Spain 2006.

\bibitem{S}
W.A.\ Strauss,
\textit{Existence of solitary waves in higher dimensions},
Comm. Math. Phys. 55 (1977), no. 2, pp.\ 149--162.

\bibitem{SS99} C. Sulem, P-L. Sulem,
\textit{The nonlinear Schr\"odinger equation. Self-focusing and wave collapse},
Applied Mathematical Sciences, 139. Springer-Verlag, New York, 1999. xvi+350 pp.

\bibitem{T04}
T.\ Tao,
\textit{On the asymptotic behavior of large radial data for a focusing non-linear
Schr\"odinger equation}, Dyn. Partial Differ. Equ. 1 (2004), no. 1, pp.\ 1--48.

\bibitem{T06}
T.\ Tao, \textit{A (concentration-)compact attractor for
high-dimensional non-linear Schr\"odinger equations}, Dynamics of PDE 4 (2007) pp.\ 1--53.

\bibitem{Tao-book}
T.\ Tao, \textit{Nonlinear dispersive equations. Local and global analysis}.
CBMS Regional Conference Series in Mathematics, 106.
Published for the Conference Board of the Mathematical Sciences, Washington, DC;
by the American Mathematical Society, Providence, RI, 2006. xvi+373 pp. ISBN: 0-8218-4143-2

\bibitem{V}
M. Vilela, \textit{Regularity of solutions to the free Schr\"odinger
equation with radial initial data},  Illinois J. Math. 45  (2001),
no. 2, pp. 361--370.

\bibitem{W83} M.\ Weinstein,
\textit{Nonlinear Schr\"odinger equations and sharp interpolation
estimates}, Comm. Math. Phys. 87 (1982/83), no. 4, pp.\ 567--576.

\bibitem{27}
V. E. Zakharov,
\textit{Collapse of Langmuir waves}, Zh. Eksp. Teor. Fiz. 62, (1972),
pp.\ 1745-1751,
(in Russian); Sov. Phys. JETP, 35 (1972), pp.\ 908-914 (English).

\end{thebibliography}
\end{document}